\documentclass[a4paper, 12pt]{article}
\usepackage{helvet}
\usepackage{amsthm}
\newtheorem{thm}{Theorem}
\newtheorem{cor}[thm]{Corollary}

\newtheorem{lemma}[thm]{Lemma}
\usepackage{amssymb,amsmath}

\newtheorem*{theorem*}{Theorem}

\DeclareMathOperator{\Q}{\mathbb{Q}}
\DeclareMathOperator{\F}{\mathbb{F}}

\begin{document}
\baselineskip=16.3pt
\parskip=14pt

\begin{center}
\section*{On the realization of subgroups of $PGL(2,F)$, and their automorphism groups, as 
Galois groups over function fields}

{\large 
Rod Gow and Gary McGuire\footnote{email {\tt gary.mcguire@ucd.ie}}
 \\ { \ }\\
School of Mathematics and Statistics\\
University College Dublin\\
Ireland}
\end{center}

 \subsection*{Abstract}
 
 Let $F$ be any field.
 We give a short and elementary proof 
 that any finite subgroup $G$ of $PGL(2,F)$ occurs as a Galois group over the function field $F(x)$.
 We also develop  a theory of descent  to subfields of $F$. This enables us to realize the automorphism groups
 of finite subgroups of $PGL(2,F)$ as Galois groups.
 \subsection*{Keywords}
 
 Galois group, linear fractional transformation

 \newpage

\section{Introduction}
  
 Given a field $F$ and a finite group $G$,
the general inverse Galois problem is to determine 
whether there exists a Galois extension of $F$ having $G$ as its Galois group.
The solution is known for certain base fields, 
such as $\mathbb{C}(x)$ where every finite group occurs as the
Galois group of some extension.
Both existence proofs and constructive proofs with explicit polynomials  appear in the literature.
The terminology
`the inverse Galois problem' usually refers to the case that the base field is $\Q$, where the 
problem is still open.
In this paper our starting point is the following theorem.

\begin{theorem*} 
 Let $F$ be any field, and let $F(x)$ be the field of rational functions in $x$ over $F$.
 Let $G$ be any finite subgroup of $PGL(2,F)$.
 Let $F(x)^G$ be the field of $G$-invariant rational functions.
 Let $f(x)/g(x)$ be a generator of $F(x)^G$.
 Then 
 $f(T)-xg(T) \in F(x)[T]$  is irreducible of degree $|G|$ and has Galois group $G$ over $F(x)$.
\end{theorem*}

 \section{Background}
 
 In this section we present some background for the later sections.
  In this paper $\overline{F}$ will denote an algebraic closure of a field $F$.
  We use the notation $\F_q$ for a finite field with $q$ elements,
  where $q$ will always be a power of a prime number $p$.
 
 \subsection{The Projective General Linear Group $PGL(2,F)$}
 
 The general linear group $GL(2,F)$ is the group of all invertible $2$-by-$2$
 matrices with entries in $F$. The projective general linear group $PGL(2,F)$
 is the quotient of $GL(2,F)$ by the subgroup consisting of all the nonzero scalar multiples of the identity
 (its centre).
 We write $PGL(2,q)$ for $PGL(2,\F_q)$. 

 \subsection{Function Fields}
 
 Let $F$ be a field,  let $x$ be an indeterminate over $F$, and let $F(x)$ be the 
 field of rational functions in $x$. We will call an element of $F(x)$ a \emph{constant} if it belongs to the base field $F$.
 A nonconstant element $\alpha$, say, of $F(x)$ is transcendental over $F$, as proved, for example, in Proposition 7.5.5 (a) of Cox.
 Thus $F(\alpha)$ is itself a field of rational functions and hence isomorphic to $F(x)$. Furthermore, $F(x)$ has finite degree over
 $F(\alpha)$. 
 
 The full group of $F$-automorphisms of $F(x)$  is  $PGL(2,F)$,
 which acts by mapping $r(x)\in F(x)$ to $r(\frac{ax+b}{cx+d})$, where $r$ is a rational function over $F$,
 $a,b,c,d \in F$, and $ad-bc\neq 0$. 
 
 Let $G$ be a finite subgroup of $PGL(2,F)$. We may thus consider $G$ to be a finite subgroup of $F$-automorphisms of $F(x)$, as explained above.
 Let $F(x)^G$ be the subfield of $F(x)$ consisting of all those elements fixed by the action of each element of $G$. Artin's theorem
 implies that $F(x)$ is a Galois extension of $F(x)^G$ with Galois group isomorphic to $G$. Thus $F(x)$ has degree $|G|$ over $F(x)^G$.
 A well known theorem of L\"uroth states that any subfield of $F(x)$ strictly containing $F$ is a rational function field over $F$ and thus
 has the form $F(y)$, where $y$ is expressible as $a(x)/b(x)$ for relatively prime polynomials $a$ and $b$ over $F$. The degree $[F(x):F(y)]$ is equal $\max(\deg(a), \deg(b))$, as proved in Proposition 7.5.5 (c) of Cox.
 
 L\"uroth's theorem in its full generality is a nontrivial result. However, in the case of the subfield $F(x)^G$, L\"uroth's
 theorem may be proved in a more elementary and constructive  manner, which also helps us to understand aspects of Galois theory
 relating to this situation, as we now explain.
 
 We consider $G$ acting on $F(x)$ as above, with $s(x)$ denoting the image of $x$ under the element $s$ of $G$. Let $T$ be another indeterminate over $F(x)$.  We define the \emph{orbit polynomial} of $G$ to be
 \[
 P_G(T)=\prod_{s\in G} \bigl(T-s(x) \bigr).
 \]
 Such polynomials are much used in Galois theory. Examples are found in earlier works, such as Section 1.1 of \cite{JLY}.
 Further information and examples can also be found in \cite{GM}.
 
 The roots of $P_G(T)$ are the elements $s(x)$, $s\in G$. These are all different. For, as $x$ is transcendental over $F$, no nonidentity
 element of $G$ can fix $x$. 
 
 By its construction, we see that the coefficients of $P_G(T)$ are fixed elementwise by each element of $G$ and thus lie in $F(x)^G$.
 Some coefficients may lie in $F$ but we claim that there are coefficients that do not lie in $F$. For $x$ is a root of $P_G(T)$ and if
 $P_G(T)$ lies in $F[T]$, this implies that $x$ is algebraic over $F$, which is not the case.
 
 Let $\Phi(x)=f(x)/g(x)$ be a nonconstant coefficient in $P_G(T)$, where $f$ and $g$ are relatively prime polynomials over $F$. Then it is shown in the proof of the special case of L\"uroth's theorem on p.23 of \cite{JLY} that $F(x)^G=F(\Phi(x))$. We call $\Phi(x)$ a \emph{generator} for $F(x)^G$ but of course it is not unique. By Proposition 7.5.5 (c), $|G|=\max(\deg(f), \deg(g))$, as we noted above. From the construction of the orbit polynomial, it is not difficult to prove that $\deg f\geq \deg g$, and thus we may take it that
 the degree of the numerator $f$ of the generator equals $|G|$. We may also arrange it so that $f$ is monic.
 
 We want to obtain a different representation for the generator $\Phi(x)$. Consider the polynomial
 \[
 f(T)-\Phi(x)g(T),
 \]
 which we view as a polynomial in $F(\Phi(x))[T]=F(x)^G[T]$. 
 
 \begin{thm} \label{orbit_polynomial_representation}
 Let $\Phi(x)=f(x)/g(x)$ be a generator of the fixed field $F(x)^G$, where $f$ is monic and has degree $|G|$. Then we have
 \[
 P_G(T)=f(T)-\Phi(x) g(T).
 \]
 \end{thm}

 \begin{proof}
 By Proposition 7.5.5 (b) of Cox, the polynomial $f(T)-\Phi(x) g(T)$ is irreducible in $F(\Phi(x))[T]$. Now since $x$ is a root of $f(T)-\Phi(x) g(T)$, and this polynomial is monic in $T$, $f(T)-\Phi(x) g(T)$ is the minimal polynomial of $x$ over $F(\Phi(x))$. Likewise,
 $x$ is a root of $P_G(T)$, whose coefficients also lie in $F(\Phi(x))$. It follows that $f(T)-\Phi(x) g(T)$ divides $P_G(T)$.
 Since both polynomials are monic and have degree $|G|$, they are identical.
 \end{proof}
 
 \begin{cor} \label{the_coefficients_of_orbit_polynomial}
 The coefficients of the orbit polynomial $P_G(T)$ are all expressible as $a+b\Phi(x)$, for suitable $a$ and $b$ in $F$.
 \end{cor}
 
 We remark that this corollary may be proved directly by examining how the coefficients of $P_G(T)$ are formed using products of elements
 of the form $s(x)$, as $s$ runs over $G$.

\section{Main Theorem}\label{sect_main}

  This section presents a proof of our main working theorem.
 
\begin{thm}\label{main1}
 Let $F$ be any field, and let $F(x)$ be the field of rational functions in $x$ over $F$.
 Let $G$ be any finite subgroup of $PGL(2,F)$. Let $\Phi(x)=f(x)/g(x)$ be a generator of the field $F(x)^G$ of fixed points of the
 $G$-action. 
 Then $G$ is the Galois group of
 $P_G(T)=f(T)-\Phi(x)g(T)$  over $F(\Phi(x))$, and $F(x)$ is the splitting field of $P_G(T)$ over $F(\Phi(x))$.
\end{thm}

\begin{proof}
The arguments are in some respects purely formal. We know that $F(x)$ is a Galois extension of $F(\Phi(x))$ with Galois group $G$.
Since we proved that $P_G(T)$ is the minimal polynomial of $x$ over $F(\Phi(x))$, the theorem follows.
\end{proof}

As we mentioned, the proof offered is largely a summary of the information we have compiled. However, we would like to recast the theorem
in what appears to be a less obvious form, although it is a little more than a translation of the theorem into a slightly different language.

\begin{thm} \label{version_over_F(x)}
Assume the notation of Theorem \ref{main1}. Let $z$ be an indeterminate over $F$. Then the polynomial
\[
f(T)-zg(T)
\]
in $F(z)[T]$ is irreducible over $F(z)$ and has Galois group $G$. Its splitting field over $F(z)$ is a rational function field
in one variable over $F$ and is thus a regular extension of $F$.
\end{thm}

\begin{proof}
Since $x$ is transcendental over $F$, and $\Phi(x)$ is not a constant, $\Phi(x)$ is transcendental over $F$, by Proposition 7.5.5 (a) of Cox.
We can thus define a field  isomorphism between $F(\Phi(x))$ and $F(z)$. This induces a ring isomorphism between $F(\Phi(x))[T]$
and $F(z)[T]$. Under this isomorphism the polynomial $f(T)-\Phi(x)g(T)$ is mapped to $f(T)-zg(T)$ and thus $f(T)-zg(T)$ is irreducible
in $F(z)[T]$. The isomorphism between $F(\Phi(x))$ and $F(z)$ then induces an isomorphism of Galois groups and splitting fields
over the stated ground fields.

It follows that the splitting field $L$, say, of $f(T)-zg(T)$ over $F(z)$ is a rational function field in one variable over $F$. We can now say that $L$ is a regular extension of $F$. For the only elements in a rational function field in one variable over $F$ that are algebraic over
$F$ are precisely the constants, again by Proposition 7.5.5 (a) of Cox.
\end{proof}

\begin{cor} \label{algebraically_closed}
Assume the notation of Theorem \ref{main1}. Let $\overline{F}$ denote the algebraic closure of $F$. Then considered as a polynomial
in $\overline{F}(z)[T]$, 
\[
f(T)-zg(T)
\]
is irreducible with Galois group $G$.
\end{cor}

\begin{proof}
We outline the details. $G$ acts as a finite group of automorphisms of the rational function field $\overline{F}(x)$ and it is clear
that the fixed field is $\overline{F}(\Phi(x))$. The rest of the proof is then a repetition of the ideas used to prove Theorems
\ref{main1} and \ref{version_over_F(x)}. 
\end{proof}



We have seen that the Galois extensions of $F(x)$ constructed by means of finite subgroups of $PGL(2,F)$ are themselves
rational functions fields of one variable over $F$. In fact, such extensions are characterized by this property, as the following simple
observation points out.

\begin{thm} \label{rational_function_field_extension}
 Let $L$ be a Galois extension of $F(x)$ of finite degree with Galois group $G$. Suppose that $L=F(z)$ is a rational function
 field in one variable $z$ over $F$. Then $G$ is isomorphic to a subgroup of $PGL(2,F)$ and $L$ is the splitting field over $F(x)$
 of the polynomial $f(T)-xg(T)$ in $F(x)[T]$, where $x=f(z)/g(z)$.
 \end{thm}
 
 \begin{proof}
 $G$ is a group of $F$-automorphisms of the rational function field $F(z)$ and hence is a subgroup of $PGL(2,F)$, as we explained
 at the beginning of this section. The rest follows from the theory we have described in this section.
 \end{proof}


We now present a few examples.

{\bf Example 1:} 
Let $F$ be any   field of characteristic not equal to 2 in which $-1$ is a square.
 For example, $F$ can be any algebraic number field containing $\Q(i)$, where $i^2=-1$,
 or any finite field $\F_q$, where $q \equiv 1 \pmod4$.
 Then $PGL(2,F)$ contains a subgroup isomorphic to $A_4$.
 One  presentation of $A_4$ given by Felix Klein \cite{K}
 occurs when the elements  are the linear fractional transformations that take $x$  to
 \[
 \biggl\{ \pm x, \pm \frac{1}{x}, \pm \frac{i(x+1)}{x-1}, \pm \frac{i(x-1)}{x+1}, 
 \pm \frac{x+i}{x-i}, \pm \frac{x-i}{x+i} \biggr\}.
 \]
 The orbit polynomial for this group is easily calculated using a computer algebra
 package such as Magma \cite{Bosma}, or by hand in a small example such as this.
A generator of the invariant functions $F(x)^{A_4}$ is then found to be $f(x)/g(x)$ where
$f(x)=x^{12} - 33x^8 - 33x^4+1$ and $g(x)=x^{10} - 2x^6 + x^2$.
By Theorem \ref{main1} 
\[
f_1(T)=T^{12} - 33T^8 - 33T^4+1-x(T^{10} - 2T^6 + T^2)
\]
is irreducible over $F(x)$ and has Galois group $A_4$.


{\bf Example 2:} 
There is another  presentation of $A_4$ as a subgroup of $PGL(2,\mathbb{C})$.
Here  $A_4$ is generated by the elements
$s(x)=\frac{x+2}{x-1}$ and $r(x)=\omega x$ where $\omega$ is a 
primitive cube root of unity.
A generator for the invariant functions is easily calculated to be
$(x^{12} + 264x^6 + 440x^3 + 24)/(x^3-1)^3$.
This remains valid over any field $F$ of characteristic not equal to 3 that contains
a primitive cube root of unity.
By Theorem \ref{main1} 
\[
f_2(T)=T^{12} + 264T^6 + 440T^3 + 24-x(T^9-3T^6+3T^3-1)
\]
is irreducible over $F(x)$ and has Galois group $A_4$.

The reason for giving two versions of $A_4$ will become apparent in the next section.

{\bf Example 3:}
Let $F$ be any field of characteristic different from 2 in which 2 is a square, and let $c=1+\sqrt{2}$.
  Then $s(x)=(-cx-1)/(x-c)$ has order 8 in $PGL(2,F)$.  
 Let $C_8$ be the cyclic subgroup generated by $s$.
 By calculating the orbit polynomial for $C_8$  we find that 
 a generator for the invariant functions is 
$ \frac{x^8 - 28x^6 + 70x^4 - 28x^2 + 1}{x^7 - 7x^5 + 7x^3 - x}$.
By Theorem  \ref{main1} 
\[
f_3(T)=T^{8} -28T^6 + 70T^4 - 28T^2 + 1 -x(T^7 - 7T^5 + 7T^3 - T)
\]
is irreducible over $F(x)$ and has Galois group $C_8$.
This is the polynomial found by Shen \cite{S} by a different method.
We obtain a  polynomial for $C_8$ over $\Q(\sqrt{2})$ for example.

{\bf Example 4:}
It is well known that the only finite subgroups of $PGL(2,\Q)$ are $C_n$ and $D_n$, 
for $n=1,2,3,4,6$.
Therefore these groups are realized over $\Q(x)$ by Theorem \ref{main1} (this is not a new result).
For example, $s(x)=\frac{-x-1}{x-1}$ generates a $C_4$.
A nonconstant coefficient in the orbit polynomial is $(x^4-6x^2+1)/(x^3-x)$, so
$T^4-6T^2+1-x(T^3-T)$ has Galois group $C_4$ over $\Q(x)$.
The same example is given  in chapter 1 of Serre \cite{S}.

As another example, the transformations $s(x)=\frac{-x-1}{x-1}$ and $t(x)=-x$ generate the dihedral group of order 8.
A nonconstant coefficient in the orbit polynomial is $(x^8+14x^4+1)/x^2(x^2-1)^2$ so the polynomial
$T^8 + 14T^4 + 1-xT^2(T^2 - 1)^2$ has Galois group $D_4$ over $\Q(x)$.

\section{Conjugacy classes and fixed points on the projective line}

 \noindent Let $G$ be a finite subgroup of $PGL(2,F)$, where $F$ is an arbitrary field. $G$ acts as a group of permutations of the
 projective line $\mathbb{P}^1(F)$. It is convenient to extend this to an action on the projective line $\mathbb{P}^1(\overline{F})$,
 where $\overline{F}$ is the algebraic closure of $F$. 
 
 
 
 We establish a number of basic principles about the existence of fixed points
 on the projective line, which are well known, but as the proofs are not complicated, we provide some
 details. Our approach partly involves looking for preferred forms, up to conjugacy, for the elements of $PGL(2,F)$.
 
 We first quote a result proved in \cite{GM}, Lemma 2.
 
 \begin{lemma} \label{cyclic_stabilizer}
Let $G$ be a finite subgroup of $PGL(2,F)$, where $F$ is any field.
Suppose that $G$ fixes a point in its action on the projective line $\mathbb{P}^1(F)$. 
\begin{enumerate}
\item If $|G|$ is not divisible by the characteristic of $F$, then $G$ is cyclic. 
\item If $|G|$ is divisible by $p=char(F)$, then 
$G$ has an elementary abelian normal $p$-subgroup with cyclic quotient of order coprime
to $p$.
\end{enumerate}
\end{lemma}






Let $\pi:GL(2,F)\to PGL(2,F)$ be the projection homomorphism. Let $F^\times$ denote the multiplicative group of nonzero elements of $F$
and let $(F^\times)^2$ denote the subgroup of squares. Let $F^{(2)}$ denote the quotient group $F^\times/(F^\times)^2$. We call
$F^{(2)}$ the group of square classes. It is an elementary abelian 2-group.

The determinant map $\det$ gives an epimorphism from $GL(2,F)$ onto $F^\times$. It then induces a homomorphism from $PGL(2,F)$
onto $F^{(2)}$ in the following way. Let $s$ be an element of $PGL(2,F)$ and write $s=\pi(g)$ for some $g\in GL(2,F)$. Define
$
\det^*: PGL(2,F)\to F^{(2)}$ by means of $\det^*(s)=(\det g)(F^\times)^2$. It is straightforward to see that $\det^*$ is a well defined
homomorphism. Its kernel is the normal subgroup $PSL(2,F)$, the projective special linear group.

 \begin{lemma} \label{pull_back}
 Let $F$ be a field and let $s$ be an element of $PGL(2,F)$ of odd finite
 order $r$. Then there exists an element $h$ of the same order in $GL(2,F)$ such that $\pi(h)=s$.
 \end{lemma}
 
 \begin{proof}
 Let $g$ be any element of $GL(2,F)$ that satisfies $\pi(g)=s$. Then $\pi(g^r)=1$ and it follows that $g^r=\lambda I$, where $\lambda$
 is a nonzero element of $F$. Let $d=\det g$. We take determinants in the equation $g^r=\lambda I$ to obtain $d^r=\lambda^2$. 
 
 Now as $s$ has odd order, and $F^{(2)}$ has exponent 2, it follows that $\det^* s=1$ and thus $d\in (F^{\times})^2$.
 We set $d=e^2$, where $e\in F$, and obtain $\lambda^2=e^{2r}$. It follows that $\lambda=\pm e^r$ and we deduce that as $r$ is odd,
 $\lambda=e^r$ or $\lambda=(-e)^r$. This argument shows that $\lambda$ is an $r$-th power in $F$. We may thus write
 $\lambda=\mu^r$, where $\mu\in F$. We then set $h=\mu^{-1}g$ and find that $h^r=I$ and $\pi(h)=s$, as required.
 \end{proof}
 
 We remark that if $r$ is even,  the conclusion of the lemma may fail. 
 
 We move towards deriving some consequences of this result when $F$ has odd prime characteristic $p$. 
 
 \begin{lemma} \label{unipotent_action}
 Let $F$ have odd prime characteristic $p$ and let $s$ be an element of $PGL(2,F)$ of order $p$. Then $s$ fixes a unique point of
 $\mathbb{P}^1(F)$ and also of  $\mathbb{P}^1(\overline{F})$.
 \end{lemma}
 
 \begin{proof}
  By Lemma \ref{pull_back}, there is an element $h$ of order $p$ in $GL(2,F)$ with $\pi(h)=s$. Since $F$ has characteristic $p$, $(h-I)^p=0$. It follows that
 the minimal polynomial $m(X)$, say, of $h$ divides $(X-1)^p$. But as $h$ is a $2\times 2$ matrix, its minimal polynomial
 has degree at most 2. Thus, $m(X)=X-1$ or $m(X)=(X-1)^2$. Certainly, $m(X)$ cannot equal $X-1$, since this would imply that
 $h$ is the identity. It follows that $m(X)=(X-1)^2$.
 
 We may consider $h$ to define a linear transformation $T$, say, of $F^2$ by setting $Tw=hw$ for 
 $w\in F^2$. (We remark here that $F^2$ denotes a two-dimensional vector space over $F$, not to be confused
 with the group $F^{(2)}$ already introduced.) Let $e$ be a vector in $F^2$ with $Te\neq e$. Let $f=Te-e$. 
 Then we have $Te=e+f$, $Tf=f$. The vectors $e$ and $f$ are clearly a basis of $F^2$ and with respect to this basis, $T$ has matrix,
 \[
 \left(
\begin{array}
{cc}
        1&0\\
        1&1 
\end{array}
\right)
\]
Let $h'$ denote this matrix of order $p$. We have shown thus that $h$ is conjugate in $GL(2,F)$ to $h'$. 

Now $\pi(h')$ defines the fractional transformation $x\mapsto x/(1+x)$. We solve the equation $s(x)=x$ to find the fixed points
of $\pi(h')$. We obtain the equation $x=x^2+x$ and hence $x=0$ is the only solution. Thus $\pi(h')$ fixes a unique point of
$\mathbb{P}^1(F)$, and since $\pi(h)=s$ is conjugate to $\pi(h')$ in $PGL(2,F)$, the same conclusion holds for $s$. The same argument shows that this is the unique point of $\mathbb{P}^1(\overline{F})$ fixed by $s$.
 \end{proof}
 
 We can now extend this argument to finite $p$-subgroups of $PGL(2,F)$. 
 
 \begin{cor} \label{finite_p_subgroup}
 Let $F$ have odd prime characteristic $p$ and let $G$ be a nontrivial finite $p$-subgroup  of $PGL(2,F)$. Then $G$ fixes a unique point of
 $\mathbb{P}^1(F)$.
 \end{cor}
 
 \begin{proof}
 Let $s$ be a nonidentity element contained in the centre of $G$. By Lemma \ref{unipotent_action}, $s$ fixes a unique point
 of $\alpha$, say, of 
 $\mathbb{P}^1(F)$. Consider now for any element $g$ of $G$ the point $g\alpha$ of $\mathbb{P}^1(F)$. We have 
 \[
 s(g(\alpha))=g(s(\alpha))=g(\alpha),
 \]
 since $g$ centralizes $s$ and $s$ fixes $\alpha$.
 
 This argument shows that $s$ also fixes $g(\alpha)$. Since we know that $s$ has a unique fixed point, we must have $g(\alpha)=\alpha$.
 Thus $G$ fixes $\alpha$ and it is clearly the unique point fixed by $G$. 
 \end{proof}
 
 Our next objective is to extend this result to the normalizer of a $p$-subgroup.
 
 \begin{cor} \label{normalizer_action}
 Let $F$ have odd prime characteristic $p$ and let $G$ be a finite subgroup of $PGL(2,F)$ of order divisible by $p$. Let $U$ be a Sylow $p$-subgroup of $G$ and let $N$ be its normalizer in $G$. Then $N$ fixes a unique point, $\alpha$, say, of $\mathbb{P}^1(F)$ and $N$
 is the full stabilizer of $\alpha$ in $G$. The stabilizer in $G$ of any point in a $G$-orbit different from that containing $\alpha$ is a (cyclic) $p'$-subgroup.
 \end{cor}
 
 \begin{proof}
 We have shown in Corollary \ref{finite_p_subgroup} that $U$ fixes a unique point, $\alpha$, say,  of $\mathbb{P}^1(F)$. Now let
 $g$ be any element of $N$ and $s$ any element of $U$. We have $g^{-1}sg(\alpha)=\alpha$, since $g^{-1}sg$ in in $U$. Thus
 \[
 s(g(\alpha))=g(\alpha),
 \]
 which shows that $s$ fixes $g(\alpha)$. By uniqueness of $\alpha$, $g(\alpha)=\alpha$, and we have shown that $N$ fixes
 $\alpha$. It is of course clear that $\alpha$ is the unique fixed point of $N$.
 
 Let $H$ be the stabilizer in $G$ of $\alpha$. We know from Lemma \ref{cyclic_stabilizer} that $H$ has a normal Sylow $p$-subgroup
 and since $N$ is necessarily a subgroup of $H$,  $U$ is also a Sylow $p$-subgroup of $H$. It follows that $U$ is normal in $H$. Thus $H$ is contained in the normalizer of $U$ and hence equals $N$. 
 
 Suppose that $\beta$ is a point of $\mathbb{P}^1(F)$ not in the same $G$-orbit as $\alpha$. Let $T$ be the stabilizer in $G$ of $\beta$.
 Suppose if possible that $p$ divides $|T|$. Let $t$ be an element of order $p$ in $T$. Sylow's theorem shows that a conjugate
 of $t$ lies in $U$ and hence this conjugate fixes $\alpha$. Then $t$ fixes some point in the $G$-orbit of $\alpha$. But $t$ also
 fixes $\beta$, which by assumption lies in a different $G$-orbit from $\alpha$. Thus $t$ has at least two fixed points,
 which we know is a contradiction. Consequently, $T$ is a $p'$-group and is hence cyclic by Lemma \ref{cyclic_stabilizer}.
 \end{proof}
 
 It will be convenient now for us to establish representative elements that describe the conjugacy classes of $PGL(2,F)$ and make calculations
 easier.
 
 \begin{lemma} \label{canonical_form}
 Let $F$ be an arbitrary field and let $A$ be an element of $GL(2,F)$ that is not a scalar multiple of the identity.
 Then $A$ is conjugate in $GL(2,F)$ to a matrix of the form
 \[
 \left(
\begin{array}
{cc}
        0&a\\
        1&b 
\end{array}
\right),
\]
 where $a$ and $b$ are elements of $F$, with $a\neq 0$.
 \end{lemma}
 
 \begin{proof}
 As in the proof of Lemma \ref{unipotent_action}, left multiplication by $A$ defines a linear transformation, $T$, say, of $F^2$. Since
 $T$ is not multiplication by a fixed scalar, there is some $u$ in $F^2$ such that $u$ and $Tu$ are linearly independent. Then $u$ and $Tu$
 are a basis of $F^2$ and with respect to this basis, the matrix of $T$ has the form above ($a$ is nonzero, as $T$ is invertible).
Thus $A$ is conjugate to this given matrix.
 \end{proof}
 
 \begin{cor} \label{nice_form}
 Let $s$ be a nonidentity element of $PGL(2,F)$. Then $s$ is conjugate in $PGL(2,F)$ to the fractional transformation
 \[
 x\mapsto 1/(cx+d),
 \]
 where $c\neq 0$, $d$  are  elements of $F$.
 \end{cor}
 
 This conjugacy result enables us to describe the involutions (elements of order two) in $PGL(2,F)$, up to conjugacy.
 
 \begin{thm} \label{involutions}
 An involution $s$ in $PGL(2,F)$ is conjugate in $PGL(2,F)$ to the fractional transformation
 \[
 x\mapsto a/x,
 \]
 where $a$ is a nonzero element of $F$. The element $s$ fixes a point of $\mathbb{P}^1(F)$ if and only if
 $a$ is a square in $F$. 
 \end{thm}
 
 \begin{proof}
 Up to conjugacy in $PGL(2,F)$, we may assume that $s$ is the fractional transformation
 \[
 s(x)=1/(cx+d),
 \]
 in accordance with Corollary \ref{nice_form}. Now we have
 \[
 s^2(x)=s(s(x))=(cx+d)/(cdx+c+d^2)=x.
 \]
 Since $x$ is transcendental over $F$, this requires that $cd=0$ and since $c\neq 0$, $d$ must be 0. We obtain the required form
 for $s$ when we set $a=c^{-1}$.
 
 Clearly, $s$ fixes a point of $\mathbb{P}^1(F)$ if and only if any conjugate of it fixes a point. The fractional
 transformation $x\mapsto a/x$ fixes a point if and only $a$ is a square in $F$, and we obtain the desired result.
 \end{proof}
 

 \begin{cor} \label{involutions_in_PSL}
 Let $F$ be a field of characteristic $2$. Suppose that $s$ is an involution in $PGL(2,F)$. Then $s$ fixes a point of $\mathbb{P}^1(F)$ if and only if $s$ is in $PSL(2,F)$.
 \end{cor}
 
 \begin{proof}
 It suffices to prove the result for involutions $s$ of the form $x\mapsto a/x$ where $a\in F^\times$. As we remarked above, this
 involution fixes a point of $\mathbb{P}^1(F)$ precisely when $a\in (F^\times)^2$. But $\det^* s=a(F^\times)^2$ and the result is clear.
 \end{proof}

The previous result concerning involutions enables us to obtain an analogue of Corollary \ref{normalizer_action} in characteristic 2.

\begin{lemma} \label{normalizer_in_characteristic_2}
Let $F$ be a field of characteristic $2$ and let $G$ be a finite subgroup of $PSL(2,F)$ of even order. Then
$G$ has an orbit on $\mathbb{P}^1(F)$ in which the stabilizer of a point is the normalizer of a Sylow $2$-subgroup. The stabilizer
of a point in any other orbit is cyclic of odd order.
\end{lemma}

\begin{proof}
Let $U$ be a Sylow $2$-subgroup of $G$ and let $t$ be an involution in the centre of $U$. As $t$ belongs to $PSL(2,F)$, Corollary
\ref{involutions_in_PSL} implies that $t$ fixes a point of $\mathbb{P}^1(F)$ and previous arguments show that as we are working
in characteristic 2, this point is unique. Again, uniqueness implies that $U$ fixes this point, and then the normalizer of $U$ also
fixes it. Finally, the proof of Corollary \ref{normalizer_action} implies that the normalizer is the stabilizer of this point (the assumption
that $p$ is odd in Corollary \ref{normalizer_action} is only needed to ensure that a fixed point exists, and we know that one exists
in our special circumstances).
\end{proof}

 We finish this section by considering the general question of determining fixed points of elements of finite order when the group $PGL(2,F)$ 
 acts on the projective line over $\overline{F}$, the algebraic closure of $F$. 
 
 \begin{thm} \label{fixed_point_theorem}
 Let $s$ be a nonidentity element of $PGL(2,F)$ of finite order. Then $s$ fixes exactly two points of $\mathbb{P}^1(\overline{F})$, unless
 $F$ has prime characteristic $p$ and $s$ has order $p$, in which case $s$ fixes exactly one point.
 \end{thm}
 
 \begin{proof}
 Corollary \ref{nice_form} shows that, up to conjugacy in $PGL(2,F)$, we may assume that $s$ is the fractional transformation
 \[
 s(x)=1/(cx+d),
 \]
 where $c$ and $d$ are elements of $F$, $c\neq 0$. To find fixed points $\alpha$ of $s$, we solve $s(\alpha)=\alpha$, which leads
 to the quadratic equation $c\alpha^2+d\alpha-1=0$.
 
 Suppose first that $F$ has characteristic 2. Then the quadratic above has two different roots, and hence $s$ has exactly two fixed points,
 if and only if $d\neq 0$. But if $d=0$, $s$ is an involution and we know that in this case, $s$ has just one fixed point.
 
 We may thus restrict our attention to when $F$ has characteristic different from 2. The quadratic  $c\alpha^2+d\alpha-1=0$ has two
 different roots provided $d^2+4c\neq 0$. We therefore consider what happens when $d^2+4c=0$. We set $e=d/2$, so that $-e^2=c$.
 In $GL(2,F)$, we set 
 \[
A=
 \left(
\begin{array}
{cc}
        1&0\\
        e&1
\end{array}
\right), \quad B=
 \left(
\begin{array}
{cc}
        e&1\\
        0&e
\end{array}
\right),\quad C=
\left(
\begin{array}
{cc}
        0&1\\
        c&d
\end{array}
\right)=\left(
\begin{array}
{cc}
        0&1\\
        -e^2&2e
\end{array}
\right).
\]
Then we calculate that $AB=CA$, and this implies that $C$ is conjugate to $B$, since $A$ is invertible. Translating this conjugacy
relation to $PGL(2,F)$, we see that $s$ is conjugate to $t$ defined by
\[
t(x)=(ex+1)/ex=x+e^{-1}.
\]

Now if $F$ has characteristic 0, $t$ has infinite order, and this is impossible, since $s$ is assumed to have finite order. If $F$ has characteristic $p$, $t$ has prime order $p$ and hence $s$ also has order $p$. We know then by Lemma \ref{unipotent_action} that
$s$ has a unique fixed point.
 \end{proof}
 
 \section{On the existence of short $G$-orbits}
 
 \noindent We assume as before that $G$ is a finite subgroup of $PGL(2,F)$. Then, as we studied in the previous section,
 $G$ acts as a permutation group of the projective lines over $F$ and $\overline{F}$. An orbit of $G$ acting on these varieties  in which the stabilizer of a point
 is the identity subgroup is called a regular or long orbit. The size of such an orbit is $|G|$. An orbit of size smaller than $|G|$ is called
 a short or nonregular orbit. The stabilizer of a point in a short orbit is larger than the identity subgroup. It seems
 to be the case that short orbits tell us more about the group and its action, perhaps because they are related to nontrivial subgroups
 and yield arithmetic information.
 
 It is a remarkable fact that $G$ has at most three short orbits on $\mathbb{P}^1(\overline{F})$. Furthermore, we know
 much about the stabilizer of a point in a short orbit. Versions of this fact can be traced back to an ingenious counting argument
 used by Felix Klein to describe the finite subgroups of the real special orthogonal group $SO(3,\mathbb{R})$, \cite{K}. We reprove
 this short orbit theorem, with special emphasis on the modular case when $F$ has prime characteristic $p$ and $p$ divides $|G|$.
 This type of analysis occurs naturally when considering the automorphism group of an algebraic curve (we are simply looking at a genus zero
 curve). 
 
 The question of when a finite group acts on a projective variety (over an algebraically closed field)
 with only finitely many short orbits seems to be one worthy of attention, but we merely raise the question here,
 and make no attempt to discuss examples or general principles.

 \begin{lemma} \label{three_non_regular_orbits}
Let $G$ be a finite subgroup of $PGL(2,F)$. Then in its action on the projective line $\mathbb{P}^1(\overline{F})$, $G$ has at most three short orbits. Moreover, if $G$ has exactly three short orbits and if $|G|$ is not divisible by the characteristic of $F$, either $G$ has a normal cyclic subgroup of index $2$ or $|G|$ is one of $12$, $24$ or $60$. 
\end{lemma}

\begin{proof}
Suppose if possible that $G$ has at least four short orbits and let $\Omega_i$, $1\leq i\leq 4$, be four such orbits. Let $\Omega$ be the union of these $\Omega_i$. Then $\Omega$ is a finite subset of $\mathbb{P}^1(\overline{F})$ on which $G$ acts with exactly four orbits. Let $G_i$ be the stabilizer of a given point in $\Omega_i$, $1\leq i\leq 4$. Note that we have $|G_i|>1$ for all $i$, since the $\Omega_i$ are not regular. 

Let $g$ be any element of $G$ and let $\chi(g)$ denote the number of fixed points of $g$ in its action on $\Omega$. We have $0\leq \chi(g)\leq 2$ if $g$ is not the identity, since such a $g$ fixes at most two points of $\mathbb{P}^1(\overline{F})$ by
Theorem \ref{fixed_point_theorem}. The well known orbit counting lemma implies that 
\[
4|G|=\sum_{g\in G} \chi(g)\leq |\Omega|+2(|G|-1).
\]
Thus since 
\[
|\Omega|=|\Omega_1|+\cdots+|\Omega_4|=\frac{|G|}{|G_1|}+\cdots+\frac{|G|}{|G_4|},
\]
we obtain the inequality 
\[
2\leq \frac{1}{|G_1|}+\cdots+\frac{1}{|G_4|}-\frac{2}{|G|}.
\]
This is clearly impossible, since each fraction $1/|G_i|$ is at most 1/2.
We deduce that $G$ has at most three nonregular orbits on $\mathbb{P}^1(F)$. 

We now consider the case that $G$ has exactly three non-regular orbits, $\Omega_i$, $1\leq i\leq 3$, with one-point stabilizers $G_i$, where we choose notation so that $|G_1|\leq \cdots \leq |G_3|$.
As before, we let $\chi(g)$ denote the number of points of $\Omega$ (the union of the $\Omega_i$) fixed by $g\in G$. The hypothesis
on the field $F$ and $|G|$ implies that $\chi(g)=2$ for all nonidentity $g$ in $G$.

The orbit counting lemma now yields that
\[
3|G|=|\Omega|+2(|G|-1)
\]
and hence
\[
1=\frac{1}{|G_1|}+\cdots+\frac{1}{|G_3|}-\frac{2}{|G|}.
\]
Given the numbering that $|G_1|\leq \cdots \leq |G_3|$, we cannot have $|G_1|\geq 3$ and thus $|G_1|=2$. This leads to
\[
\frac{1}{2}=\frac{1}{|G_2|}+\frac{1}{|G_3|}-\frac{2}{|G|}.
\]

It is clear from this inequality that $|G_2|=2$ or 3. If $|G_2|=2$, we obtain $2|G_3|=|G|$.
Thus $G_3$ has index 2 in $G$ and is hence normal. Furthermore, as $G_3$ fixes a point in
$\Omega_3$, it is cyclic by Lemma \ref{cyclic_stabilizer}. Thus $G$ has a normal cyclic subgroup of index 2 in this case.

Consider next the other possibility that $|G_2|=3$. Then we obtain
\[
\frac{1}{6}=\frac{1}{|G_3|}-\frac{2}{|G|}.
\]
From this equality, it is clear that $3\leq |G_3|\leq 5$ and the three feasible values for
$|G_3|$ give the three outcomes $|G|=12, 24$ or 60.

\end{proof}


We remark that the group theoretic information contained in our proof enables us to identify the 
groups concerned, using the group  action on the cosets  of the stabilizer subgroups.
In the cases $|G|=12, 24$ or 60, it can be shown that $G$ is isomorphic to $A_4$, $S_4$ or $A_5$
respectively. We omit the details.

We move towards the modular case when the characteristic of $F$ divides $|G|$. Parts of the lemma below appeared as Lemma 3 of \cite{GM}, but with a slightly more specific hypothesis on the field $F$.

 \begin{lemma} \label{Sylow}
Let $F$ be a field of prime characteristic $p$ and let 
$G$ be a finite subgroup of $PGL(2,F)$, where $p$ divides $|G|$. Let $Q$ be a Sylow $p$-subgroup of $G$
and let $N_G(Q)$ denote its normalizer in $G$. 
Then the following are true:

\begin{enumerate}

\item $Q$ is elementary abelian.

\item The centralizer in $G$ of any non-identity element of $Q$ is $Q$.

\item Any two different Sylow $p$-subgroups of $G$ have trivial intersection. 

\item Let $M$ denote the number of elements of order $p$ in $G$. Then 
\[
M=|G:N_G(Q)|(|Q|-1).
\]

\end{enumerate}
\end{lemma}

\begin{proof} 

1.  $Q$ fixes a point of $\mathbb{P}^1(\overline{F})$ by Corollary \ref{finite_p_subgroup} if $p$ is odd. 
The same proof works when $p=2$ since any element of order two fixes a unique point of $\mathbb{P}^1(\overline{F})$ .
 It now follows from Lemma  \ref{cyclic_stabilizer}
that $Q$ is elementary abelian, since it is a subgroup of the stabilizer of a point.

2. Let $s$ be a nonidentity element of $Q$. It is immediate that $s$ has order $p$, since $Q$ has exponent $p$.
Let $t$ be a nonidentity of $G$ that commutes with $s$. Let $\alpha$ be the unique element of $\mathbb{P}^1(\overline{F})$ that
is fixed by $s$. It is clear that $s$ fixes $t(\alpha)$, since $s$ commutes with $t$. Hence, by uniqueness of $\alpha$,
$t$ fixes $\alpha$. 

Suppose if possible that $t$ fixes another element $\beta$ of $\mathbb{P}^1(\overline{F})$ different from $\alpha$. Then
$s(\beta)$ is also fixed by $t$. Since by Theorem \ref{fixed_point_theorem}, $t$ fixes at most two points, which must be $\alpha$ and $\beta$,
$s(\beta)$ is either $\alpha$ or $\beta$. If $s(\beta)=\alpha$, then $\beta=s^{-1}(\alpha)=\alpha$, and we have a contradiction.
If $s(\beta)=\beta$, $s$ fixes two different points of of $\mathbb{P}^1(\overline{F})$ and this is another contradiction.

It follows that $t$ fixes a unique point of of $\mathbb{P}^1(\overline{F})$  and hence has order $p$ by Theorem \ref{fixed_point_theorem}.
Thus the centralizer of $s$ in $G$ is a $p$-group. Since $Q$ is abelian, $Q$ is contained in the centralizer of $s$ and thus be
the centralizer, since it is a Sylow $p$-subgroup of $G$.

3. We may assume that $Q$ is not the unique Sylow $p$-subgroup of $G$, for otherwise there is nothing to prove. Let $Q_1$ be a different Sylow $p$-subgroup and let
$g$ be an element of $Q\cap Q_1$. Then both $Q$ and $Q_1$ centralize $g$ and hence $g$ must be the identity, by what we proved in the previous paragraph.

4. To count the number of elements of order $p$, we recall that any element of order $p$ is in some Sylow $p$-subgroup, and this Sylow subgroup is unique, by our argument above. 
Furthermore, all Sylow $p$-subgroups are conjugate in $G$ and their number is $|G:N_G(Q)|$. 
Given this data, our statement about the number of elements of order $p$ is immediate. 
 \end{proof}

 This data enables us to analyze the modular case.

\begin{thm} \label{two_non_regular_orbits}
Let $F$ be a field of prime characteristic $p$ and let 
$G$ be a finite subgroup of $PGL(2,F)$, where $p$ divides $|G|$.  Then in its action on  $\mathbb{P}^1(\overline{F})$, $G$ has at most two short orbits. $G$ has exactly one short orbit if and only if $G$ is a $p$-group.
\end{thm}

\begin{proof}
We employ the notation previously used. Suppose by way of contradiction that $G$ has exactly three 
short orbits, $\Omega_i$, $1\leq i\leq 3$, with one-point stabilizers $G_i$, where $|G_1|\leq |G_2|\leq |G_3|$.
Let $M$ be the number of elements of order $p$ in $G$. The orbit counting theorem gives
\[
3|G|=|\Omega|+M+2(|G|-M-1)
\]
and since $\Omega$ is the union of three $G$-orbits, we obtain
\[
|G|+M+2=|\Omega|=\frac{|G|}{|G_1}+\frac{|G|}{|G_2|}+\frac{|G|}{|G_3|}.
\]

The previous arguments show that $|G_1|=2$ and $2\leq |G_2|\leq 3$. Furthermore, if $Q>1$ is a Sylow
$p$-subgroup of $G$, Corollary \ref{normalizer_action}  shows that we may take
$N_G(Q)$ to be one of the $G_i$.

We consider the possibility that 
$|G_2|=2$ also. We claim that $N_G(Q)$ cannot be $G_1$ or $G_2$. For if $N_G(Q)$ is one of $G_1$ or $G_2$, then $p=2$, and $G_1$ and $G_2$ are both Sylow 2-subgroups of $G$, of order 2. But then $G_1$ and $G_2$ are conjugate in $G$, and hence an element of order 2 in $G$ fixes points in both $\Omega_1$ 
and $\Omega_2$, contrary to the fact that an element of order $p=2$ fixes only one point of the projective line.

Thus we can assume in these circumstances that $G_3=N_G(Q)$. Thus, substituting $|G_1|=|G_2|=2$ and $|G_3|=|N_G(Q)|$, the orbit counting theorem yields
\[
(M+2)|G_3|=|G|.
\]
When we apply Lemma \ref{Sylow}, we now have 
\[
M+2=|G:N_G(Q)|, \quad M=|G:N_G(Q)|(|Q|-1)
\]
and these two equations are clearly incompatible. 

We are left with the possibilities that $|G_2|=3$ and $3\leq |G_3|\leq 5$ (but we are no longer assuming that
$G_3=N_G(Q)$). Corresponding to the three values for $|G_3|$, there are accompanying equalities
\begin{equation}\label{eq1}
6(M+2)=|G|, \quad 12(M+2)=|G|, \quad 30(M+2)=|G|.
\end{equation}
Furthermore, since $p$ divides $|G_i|$ for some $i$, $p$ is one of 2, 3 and 5.

Suppose first that $p=5$. Then we can only have $|G_3|=5$ and $G_3=N_G(Q)$. Our formula for
$M$ yields that $M=4|G|/5$ and this is incompatible with $30(M+2)=|G|$. 

Next, consider the case that $p=3$. Since no $|G_i|$ is divisible by 9, it follows that $|Q|=3$. Moreover, if
$N_G(Q)$ is not equal to $Q$, its order is at least 6. Thus, since no $G_i$ has order greater than 5, we  deduce that $N_G(Q)=Q$ here and $M=2|G|/3$. This is incompatible with all of the equalities \eqref{eq1},
one of which must hold true. This excludes the case that $p=3$.

The possibility that $p=2$ remains to be excluded. Now $M+2$ is odd, since the number of involutions in a finite group
of even order is odd, and thus a Sylow
2-subgroup $Q$ of $G$ has order at most 4. Suppose that $|Q|=4$. Then we must have $N_G(Q)=Q$ and since both 6 and 30 are not divisible by 4, $12(M+2)=|G|$ must hold. Furthermore,
$M=3|G|/4$. This clearly leads to a contradiction. 

Finally, suppose that $|Q|=2$. We again must have $N_G(Q)=Q$ and then $M=|G|/2$. This also leads immediately to a contradiction, 
and we deduce that there is no example of a group with exactly three nonregular orbits and order divisible by $p$.

We now consider the case that $G$ has exactly one non-regular orbit, $\Omega$, say. Our previous arguments imply that
\[
|G|=|\Omega|+M+2(|G|-M-1)
\]
and thus 
\[
M=|\Omega|+|G|-2.
\]
Now we clearly have the trivial inequality $M\leq |G|-1$, with equality only if $G$ is a $p$-group, and it follows that $G$ is indeed a $p$-group and $|\Omega|=1$ (so that $\Omega$ consists of a single $G$-fixed point).
\end{proof}

\section{Galois groups of polynomials related to $\Phi(x)$}

\noindent We return in this section to an hypothesis considered at the beginning of the paper, which
we briefly recall here. $F$ is a field and $F(x)$ is the field of rational functions
in the indeterminate $x$ over $F$. $G$ is a finite subgroup of $PGL(2,F)$. $G$ acts as a group of $F$-automorphisms of $F(x)$ and $F(x)^G$
denotes the subfield of $G$-invariant elements of $F(x)$. We have $F(x)^G=F(\Phi(x))$, where $f(x)$ and $g(x)$ are relatively
prime elements of $F[x]$, and  $\deg f=|G|$. 

We have shown in Theorem \ref{version_over_F(x)} that the polynomial $f(T)-xg(T)$ is irreducible in $F(x)[T]$ and its Galois group over $F(x)$ is isomorphic to $G$.
This section is concerned with the idea of specializing the indeterminate $x$ to an element $\lambda$ of $F$ and thus
forming the polynomial $\Phi_\lambda=f(T)-\lambda g(T)$ in $F[T]$. We show that if $\alpha$ is root of $\Phi_\lambda$, $F(\alpha)$ is always a separable extension of $F$, except possibly when $F$ has characteristic 2 and $G$ has a subgroup (necessarily normal) of index 2. 
Our aim then is to investigate the Galois
group of $\Phi_\lambda$ over $F$ in the separable situation (this is the same as the Galois group of
$F(\alpha)$ over $F$) and our main finding is that the Galois group is isomorphic to a subgroup of $G$. There
are some interesting special cases that arise when $\alpha$ is in a short orbit of the $G$-action on $\mathbb{P}^1(\overline{F})$. Here,
the Galois group of $F(\alpha)$ over $F$ has order at most 2.

It may be argued that a theorem that identifies the Galois group of $\Phi_\lambda$ merely with some (unidentified) subgroup
of $G$ is too imprecise to be useful. We feel that this description cannot be improved, as we make no assumptions on the field $F$, which could say be algebraically closed, when all Galois groups trivial. When $F$ is finite, the Galois group has to be cyclic, although $G$ may for example be a nonabelian simple group such as $PSL(2,q)$.

There are fields $F$, known as Hilbertian fields, for which there are infinitely many specializations of $x$ to an element
$\lambda$ of $F$ such that the resulting polynomial $\Phi_\lambda$ is irreducible with Galois group $G$. Hilbertian fields
include algebraic number fields and finite extensions of rational function fields $\mathbb{F}_q(x)$ (so called global fields).
As far as we can see, our investigations tell us nothing about  suitable specializations that will retain the Galois group
$G$.

We begin by enunciating the basic principle for investigating $\Phi_\lambda$. 

\begin{thm} \label{multiplicity_theorem_for_the_roots}
Assume the hypotheses given at the beginning of this section. Let $\alpha$ be a root of $\Phi_\lambda$ in $\overline{F}$ (thus
$\lambda=\Phi(\alpha)$) and let $\Gamma$ be the stabilizer in $G$ of $\alpha$ when $G$ acts on $\mathbb{P}^1(\overline{F})$. Then
we have the factorization
\[
f(T)-\lambda g(T)=f(T)-\Phi(\alpha) g(T)=\prod_{s\in G/\Gamma}(T-s(\alpha))^{|\Gamma|},
\]
where $G/\Gamma$ is a set of left coset representatives of $\Gamma$ in $G$.
\end{thm}

This result appears as Lemma 9 of \cite{GM} and also as Proposition 3.8 of \cite{Bl}.
We first apply it when $\alpha$ above is in a regular $G$-orbit, this being the generic case.


\begin{thm} \label{basics_on_Galois_groups}
Let $G$ be a finite subgroup of $PGL(2,F)$ and let $\Phi(x)=f(x)/g(x)$ be a generator of $F(x)^G$. 
Given $\lambda\in F$, set $\Phi_\lambda(T)=f(T)-\lambda g(T)\in F[T]$. Let $\alpha$ be a root of $\Phi_\lambda$
in $\overline{F}$ and suppose that $\alpha$ is contained in a regular orbit of $G$ in its action on $\mathbb{P}^1(\overline{F})$. 
Then all the irreducible factors of $\Phi_\lambda$ in $F[T]$ have the same splitting field, $F(\alpha)$, and hence the same Galois
group over $F$. The Galois group is isomorphic to a subgroup of $G$ and $|G|/|H|$ equals the number of irreducible
factors of $\Phi_\lambda$.
\end{thm}

\begin{proof}
 The $G$-orbit of 
$\alpha$ consists of $|G|$ different elements of $\overline{F}$, since $\alpha$ lies in a regular $G$-orbit by assumption. The $G$-orbit
of $\alpha$ thus consists of all the roots of $\Phi_\lambda$, since $\Phi_\lambda$ has degree $|G|$. We note also that
all roots of $\Phi_\lambda$ lie in $F(\alpha)$.

Let $\phi(T)$ be an irreducible factor of $\Phi_\lambda$ in $F[T]$. The roots of $\phi$ are a subset of those of $\Phi_\lambda$, hence  they
all lie in $F(\alpha)$, and we deduce that $F(\alpha)$ is a splitting field for each irreducible factor. Let $H$ be the Galois group
of $F(\alpha)$ over $F$. This is the same as the Galois group of each irreducible factor of $\Phi_\lambda$ and of $\Phi_\lambda$ itself.

Let $\sigma$ be an element of $H$ and $s$ an element of $G$. The elements of $H$ act trivially on $F$ and hence from the definition of fractional transformations over $F$,
\[
\sigma(s(\alpha))=s(\sigma(\alpha)).
\]
Thus the actions of $H$ and $G$ on the roots of $\Phi_\lambda$ commute. It follows from Lemma 8 of \cite{GM} that $G$ transitively permutes
the orbits of $H$ acting on the roots of $\Phi_\lambda$ and if $S$ is the stabilizer in $G$ of an orbit, $S$ is isomorphic to 
$H$. This of course means that $H$ is isomorphic to a subgroup of $G$. Moreover, the different $H$-orbits on the roots of $\Phi_\lambda$ correspond to the different irreducible factors of $\Phi_\lambda$ in $F[T]$, and their number is the index of $S$ in $G$, which equals $|G|/|H|$.
\end{proof}

We note that a version of this result, stated for finite fields, appears as Theorem 24 of \cite{GM}. The proof is essentially the same.

Our next requirement is to extend Theorem \ref{basics_on_Galois_groups} to deal with roots lying in short orbits. Initially
we meet with potential separability problems but these are confined to fields of characteristic 2 and indeed to groups $G$ that have
a normal subgroup of index 2.

\begin{thm} \label{separable_extension}
Let $G$ be a finite subgroup of $PGL(2,F)$ and let $\Phi(x)=f(x)/g(x)$ be a generator of $F(x)^G$. 
Given $\lambda\in F$, set $\Phi_\lambda(T)=f(T)-\lambda g(T)\in F[T]$. Let $\alpha$ be a root of $\Phi_\lambda$
in $\overline{F}$ and suppose that $\alpha$ is contained in a short orbit of $G$ in its action on $\mathbb{P}^1(\overline{F})$. 
Suppose that $F(\alpha)$ is an inseparable extension of $F$. Then $F$ has characteristic $2$ and $G$ contains a (normal) subgroup
of index $2$.
\end{thm}

\begin{proof}
It is certainly the case that $F$ has prime characteristic $p$ if it has an inseparable extension. Furthermore, the degree of $F(\alpha)$
over $F$ must be divisible by $p$. Now each irreducible factor of $\Phi_\lambda$ has splitting field $F(\alpha)$ and hence the roots
in $F(\alpha)$ of such a factor occur with multiplicity divisible by $p$.
Let $\Gamma$ be the stabilizer of $\alpha$ in $G$. Then since $\alpha$ occurs with multiplicity
$|\Gamma|$ as a root of $\Phi_\lambda$, $p$ divides $|\Gamma|$. Now as we have seen, there is a unique orbit of $G$ acting on
$\mathbb{P}^1(\overline{F})$ whose stabilizer has order divisible by $p$. In this case, $\Gamma$ is the normalizer of a Sylow $p$-subgroup
of $G$. Now we have shown that this short orbit is defined on $\mathbb{P}^1(F)$ if $p$ is odd or if $p=2$ and $G$ has no subgroup
of index 2. This means that $\alpha$ is an element of $F$, and hence the extension is (trivially) separable. It follows that $p=2$ in the event
of an inseparable extension.
\end{proof}


\begin{thm} \label{Galois_group_of_order_two}
Let $G$ be a finite subgroup of $PGL(2,F)$ and let $\Phi(x)=f(x)/g(x)$ be a generator of $F(x)^G$. 
Given $\lambda\in F$, set $\Phi_\lambda(T)=f(T)-\lambda g(T)\in F[T]$. Let $\alpha$ be a root of $\Phi_\lambda$
in $\overline{F}$ and suppose that $F(\alpha)$ is a Galois extension of $F$ (this is guaranteed unless possibly
$F$ has characteristic $2$). Suppose that $\alpha$ is contained in a short orbit of $G$ in its action on $\mathbb{P}^1(\overline{F})$. 
Let $\Gamma$ be the stabilizer of $\alpha$ in $G$ and let $H$ be the Galois group of $F(\alpha)$ over $F$. Then there is a subgroup,
$S$, say, of $G$ containing $\Gamma$ as a normal subgroup with $S/\Gamma$ isomorphic to $H$. Furthermore, $H$ has order at most $2$.
\end{thm}

\begin{proof}
Let $\Omega$ be the set of different roots of $\Phi_\lambda$ in $\overline{F}$. We know that $G$ acts transitively on $\Omega$ with
one point stabilizer equal to $\Gamma$. Since the action of $H$ on $\Omega$ commutes with that of $G$, $\Omega$ is the union of $H$-orbits
that are transitively permuted by $G$. 

Let $\Omega'$ be the $H$-orbit that contains $\alpha$ and let $S$ be the stabilizer of $\Omega'$ in $G$. Certainly, $\Gamma$ is a subgroup
of $S$ and we claim that $\Gamma$ acts trivially on $\Omega'$. For take $\beta$ in $\Omega'$. Then we can write $\beta=\sigma(\alpha)$
for some $\sigma\in H$. Then if $t\in \Gamma$, we have
\[
t(\beta)=t(\sigma(\alpha))=\sigma(t(\alpha))=\sigma(\alpha)=\beta.
\]
This shows that $\Gamma$ acts trivially on $\Omega'$. But clearly, if $s\in S$ acts trivially on $\Omega'$, $s$ fixes $\alpha$ and
thus $s$ in in $\Gamma$. We have thus shown that $\Gamma$ is the kernel of the action of $S$ on $\Omega'$ and is thus normal
in $S$.

We therefore have an action of the quotient group $S/\Gamma$ on $\Omega'$ and the argument above shows that this action is regular. Since the actions
of $H$ and $S/\Gamma$ on $\Omega'$ commute and $S/\Gamma$ acts regularly, $H$ is isomorphic to $S/\Gamma$, by Lemma 8 of \cite{GM}.

Finally, we have shown that the stabilizer $\Gamma$ of a point in a short $G$-orbit is cyclic except if $F$ has prime characteristic $p$
and $\Gamma$ is the normalizer of a Sylow $p$-subgroup of $G$. We consider these two possibilities separately. Suppose
that $\Gamma$ is cyclic. Then we know that $\Gamma$ fixes at most two points of $\mathbb{P}^1(\overline{F})$. Let $t$ be any element
of $S$ and let $s$ be an element of $\Gamma$. Since $S$ normalizes $\Gamma$, $s$ fixes $t(\alpha)$. Thus there are at most
two points in the $S$-orbit of $\alpha$ and hence $S/\Gamma$ has order at most 2.

On the other hand, if $F$ has characteristic $p$ and $\Gamma$ is the normalizer of a Sylow $p$-subgroup of $G$, it is clear
that the normalizer is self-normalizing in $G$ and hence $\Gamma=S$, which implies the Galois group is trivial.  (We remark
that this part of the proof can be replaced by an argument already used in Theorem \ref{separable_extension}.)
\end{proof}

The following example illustrates aspects of this theorem when $F$ is a finite field.

{\bf Example:} Take $F=\mathbb{F}_q$, where $q$ is a power of a prime, and $G=PGL(2,q)$. A generator for $\mathbb{F}_q(x)^G$ is
$\Phi(x)=f(x)/g(x)$, where
\[
f(x)=(1+(x^q-x)^{q-1})^{q+1},\quad g(x)=(x^q-x)^{q(q-1)}.
\]
A root $\alpha$ of $\Phi_0=f(x)$ satisfies
$-1=(\alpha^q-\alpha)^{q-1}$. In the action of $G$ on the roots of $\Phi_0$, $\alpha$ is in a short orbit whose stabilizer
is cyclic of order $q+1$. It is easy to see that $\Phi_0$ is the product of all monic irreducible polynomials of degree 2 over
$\mathbb{F}_q$ and all corresponding Galois groups have order 2. 

We turn to consideration of exceptional behaviour exhibited by the roots of the polynomial $\Phi_\lambda$ in a nonseparable
situation, which we have already shown can only happen in characteristic 2.

{\bf Example:} Let $E$ be the function field $\mathbb{F}_2(z)$, where $z$ is transcendental over $\mathbb{F}_2$. Let $F=E(\alpha)$,
where $\alpha$ satisfies $\alpha^2+\alpha=z+1$. $F$ is a separable extension of $E$ of degree 2 and we can check that
$z$ is not a square in $F$. Consider the elements $s$ and $t$ of $PGL(2,F)$ defined by the fractional transformations
\[
s(x)=\frac{\alpha x+z}{x+\alpha+1},\quad t(x)=\frac{z}{x}.
\]
We find that $s$ has order 3, $t$ has order 2, and $t^{-1}st=s^{-1}$. Thus $s$ and $t$ generate a subgroup $G$, say, of $PGL(2,F)$ isomorphic to the symmetric group $S_3$. 

We can check that $G$ has no short orbits on $\mathbb{P}^1(F)$. This shows that Corollary \ref{normalizer_action} does not
extend to characteristic 2. The only fixed point of $t$ on $\mathbb{\overline{F}}$ is $\sqrt{z}$, and $F(\sqrt{z})$ is of course
an inseparable extension. 

Let $\Phi(x)=f(x)/g(x)$ be a generator of $F(x)^G$. Then $\sqrt{z}$ is not a root of $g$ (since it is not in the orbit of $\infty$)
and thus $\lambda=\Phi(\sqrt{z})$ is an element of $\overline{F}$. If we set $\alpha=\sqrt{z}$, $\alpha$ is a root of
$\Phi_\lambda$ and the field extension $F(\alpha)$ is inseparable over $F$.

\begin{proof}
We know from the argument of Corollary \ref{normalizer_action} that $N$ fixes a unique point of $\mathbb{P}^1(F)$. Now as $|N|$ is even,
$N$ contains an involution, $t$, say, and $t$ must therefore fix a point of $\mathbb{P}^1(F)$. Corollary \ref{involutions_in_PSL} implies
that $-1$ is a square in $F$. Suppose that $p\equiv 3\bmod 4$. Then an element $\alpha$ say in $F$ that satisfies $\alpha^2=-1$ is not
in $\mathbb{F}_p$ and thus $\mathbb{F}_p(\alpha)$ is the field of order $p^2$.
\end{proof}

 \begin{thm} \label{A_4_subgroup}
 Let $F$ be a field of characteristic $2$. Suppose that $PGL(2,F)$ contains a subgroup $G$ isomorphic to $A_4$. Then $G$ fixes a unique
 point of $\mathbb{P}^1(F)$ and $F$ contains a copy of the field $\mathbb{F}_4$. Up to conjugacy in $PGL(2,F)$, $G$ consists
 of affine transformations $x\to ax+b$, where $a$, $b$ are in $F$, $a^3=1$.
 \end{thm}
 
 \begin{proof}
 We note that $A_4$ contains no subgroup of index 2. It follows that the determinant homomorphism $\det^*$ restricted to $G$ is trivial
 since the image of $\det^*$ is a group of exponent 2. (This amounts to saying that $G$ is contained in $PSL(2,F)$.)
 
  Corollary \ref{involutions_in_PSL} implies that each involution in $G$ fixes a unique point of $\mathbb{P}^1(F)$ and hence
  the Sylow 2-subgroup of $G$ fixes the same point. Since the Sylow
 2-subgroup of $G$ is normal in $G$, $G$ fixes a unique point of $\mathbb{P}^1(F)$.
 
 We may assume that the fixed point of $G$ is $\infty$ by replacing $G$ with a conjugate subgroup in $PGL(2,F)$, if necessary. Then the elements of $G$ have the
 form $s_{a,b}$, where $s_{a,b}(x)=ax+b$, $a\in F^\times$, $b\in F$. As we showed  in the proof of Lemma \ref{cyclic_stabilizer},
 the map sending $s_{a,b}$ to $a$ is a homomorphism from $G$ into $F^\times$, whose kernel
 is the Sylow 2-subgroup of $G$. It follows that the image has order 3 in $F^\times$, and hence $F^\times$ has an element of order 3.
 Consequently, $F$ contains $\mathbb{F}_4$, which is generated over $\mathbb{F}_2$ by an element of order 3.
 \end{proof}
 
 \begin{cor} \label{A_5_subgroup}
 Let $F$ be a field of characteristic $2$. Then $PGL(2,F)$  contains a subgroup isomorphic to $A_5$ if and only if $F$ contains a copy of the field $\mathbb{F}_4$. 
 \end{cor}
 
 \begin{proof}
 Suppose that $G$ contains a subgroup isomorphic to $A_5$. Then it also contains a subgroup isomorphic to $A_4$, and the conclusion
 that $F$ contains a copy of $\mathbb{F}_4$ follows from Theorem \ref{A_4_subgroup}.
 
 Conversely, suppose that $F$ contains a copy of $\mathbb{F}_4$. Then $PGL(2,F)$ contains a subgroup isomorphic to $PGL(2, 4)$,
 which is well known to be isomorphic to $A_5$.
 \end{proof}

\section{On the nature of the function $\Phi(x)=f(x)/g(x)$}

\noindent Let $G$ be a finite subgroup of $PGL(2,F)$ and let $\Phi(x)=f(x)/g(x)$ be a generator of $F(x)^G$. We assume, as is permitted,
that $\deg f=|G|$ and $\deg g<|G|$. The denominator polynomial $g$ factors completely over $F$ and we aim in this section to describe
its factors and their multiplicities. 

Let the elements of $G$ be $s_i$, $1\leq i\leq |G|$, with
\[
s_i(x)=\frac{a_ix+b_i}{c_ix+d_i},
\]
where $a_i$, $b_i$, $c_i$ and $d_i$ are in $F$, $a_id_i-b_ic_i\neq 0$, for all $i$. We take $s_1$ to be the identity,
so that $a_1=d_1=1$, $b_1=c_1=0$.

The following result should be obvious.

\begin{lemma} \label{orbit_of_infty}
In the action of $G$ on $\mathbb{P}^1(F)$, the orbit of $\infty$ consists of $\infty$ and the different elements among $a_ic_i^{-1}$. (Note
that if $c_i=0$, $a_ic_i^{-1}$ is $\infty$ and thus can be excluded.)
\end{lemma}

\begin{cor} \label{orbit_of_infty_is_regular}
If the $G$-orbit of $\infty$ is regular, no $c_i=0$ for $i>1$.
\end{cor}

In connection with the assumption on regularity made above, we note the following. Suppose that $F$ is infinite. Then as we have proved that
as $G$ has at most three nonregular orbits on $\mathbb{P}^1(F)$, $G$ certainly has regular orbits on $\mathbb{P}^1(F)$. Since
$PGL(2,F)$ acts transitively on $\mathbb{P}^1(F)$, we can replace $G$ by a conjugate under the action of $PGL(2,F)$
for which the orbit containing $\infty$ is regular. Thus, at the expense of using a conjugate, the hypothesis may be fulfilled.
However, properties concerning the invariant functions $\Phi(x)$ may be lost. On the other hand, if $F$ is finite, this procedure
is often impossible.

We shall now work with the hypothesis that the $G$-orbit of $\infty$ is regular. Then we know that the $c_i$ for $i>1$ are all nonzero
and hence by rescaling the fractional transformations, we may suppose that the elements of $G$ are given by
\[
s_i(x)=\frac{a_ix+b_i}{x+d_i}
\]
for $i>1$.

\begin{lemma} \label{another_description_of_the_orbit_of_infty}
Suppose that the $G$-orbit of $\infty$ is regular. Then the elements $d_i$, $i>1$, above are all different, and together with
$\infty$, the elements $-d_i$ constitute the $G$-orbit containing $\infty$.
\end{lemma}

\begin{proof}
We have
$s_i(-d_i)=\infty$ for $i>1$.
Thus the $G$-orbit containing $\infty$ consists of $\infty$ and the elements $-d_i$.
It follows that the elements $d_i$ are all different for $i>1$, since the orbit has size $|G|$.
\end{proof}

Under the regularity hypothesis, we can now such that a certain term arising as a coefficient of the orbit
polynomial of $P_G(T)$ is nonzero and its denominator determined.

\begin{thm} \label{the_trace_term_of_the_orbit_polynomial}
Let $G$ be a finite subgroup of $PGL(2,F)$. and let $s_i$, $1\leq i\leq |G|$, be the elements of $G$.
Suppose that the $G$-orbit of $\infty$ is regular. Let the elements $s_i$ of $G$ be represented as fractional
transformations $s_i(x)=(a_ix+b_i)/(x+d_i)$ for $i>1$, and $s_1(x)=x$. Let
\[
P_G(T)=\prod_{s\in G} \bigl(T-s(x) \bigr)
\]
be the orbit polynomial of $G$. Then the coefficient of $T^{|G|-1}$ in $P_G(T)$ (the so-called
trace term) has the form $a(x)/b(x)$, where $a(x)$ is a polynomial of degree $|G|$ and
\[
b(x)=\prod_{i>1} \bigl(x+d_i\bigr).
\]
\end{thm}

\begin{proof}
Let $S$ be the coefficient of $T^{|G|-1}$ in $P_G(T)$. Then
\[
-S=\sum_{i\geq 1}s_i(x)=x+\sum_{i\geq 2} (a_ix+b_i)/(x+d_i).
\]
We set
\[
R(x)=\prod_{i\geq 2} \bigl(x+d_i\bigr), \quad R_i(x)=R(x)/(x+d_i), \quad i\geq 2.
\]
A straightforward manipulation shows that
\[
-SR(x)=xR(x)+\sum_{i\geq 2}(a_ix+b_i)R_i(x).
\]

We claim that no linear polynomial $x+d_j$ divides the right hand side of the equality above. For since $x+d_j$ divides $R(x)$ and also
$R_i(x)$ for $i\neq j$, if $x+d_j$ divides the right hand side, it divides $(a_jx+b_j)R_j(x)$. It is clear that $x+d_j$ does not divide
$R_j(x)$, since the $d_i$ are all different under the regularity hypothesis. It follows that $x+d_j$ is a scalar multiple
of $a_jx+b_j$, and this is clearly impossible as $s_j$ is a fractional transformation. 

Now we have
\[
-S=\frac{xR(x)+\sum_{i\geq 2}(a_ix+b_i)R_i(x)}{R(x)}
\]
and we have shown that the numerator above is relatively prime to the denominator $R(x)$. The numerator has degree $|G|$
and this proves what we wanted.
\end{proof}

\begin{cor} \label{the_roots_of_g(x)}
Under the hypothesis of Theorem \ref{the_trace_term_of_the_orbit_polynomial}, a generator $\Phi(x)$ of $F(x)^G$ is expressible as
$f(x)/g(x)$, where $f(x)$ has degree $|G|$ and the denominator $g(x)$ is the product 
\[
\prod_{i\geq 2}(x+d_i).
\]
The roots $-d_i$ of $g(x)$ together with $\infty$ are a regular $G$-orbit.
\end{cor}

\section{Automorphisms of Base Field and Descent}\label{sect_descent}

There is an interesting  point concerning rationality and fields of definition that we wish to  discuss in this section.
Consider the examples in the previous section when the base field $F$ 
is  $\Q(i)$ in example 1,  $\Q(\omega)$ in example 2, and  $\Q(\sqrt{2})$ in example 3.
In all cases, the polynomial $f_i(T)$ actually has coefficients in $\Q(x)$.
Therefore we can ask for the Galois group of $f_i(T)$ over $\Q(x)$. Analogous phenomena hold for fields such as $\mathbb{F}_q(x)$.
This is  the descent problem already described and we may be able 
to realize groups other than subgroups of $PGL(2,F)$
as Galois groups over function fields by working over subfields of  $F$.
Such ideas have been pursued before; see, for example, Section 4 of \cite{Ab2}.

We consider the following general configuration. Let $F$ be a field and let $H$ be a finite group of automorphisms of $F$. Let $E$ be the fixed field of $H$. We may extend the action of $H$ from $F$ to $F(x)$ coefficientwise.
Explicitly,   if $\sigma\in H$ and if $(\sum a_i x^i )/(\sum b_i x^i) \in F(x)$
then 
\[
\sigma \biggl( \frac{\sum a_ix^i}{\sum b_i x^i}\biggl) =
 \frac{\sum \sigma( a_i) x^i}{\sum \sigma(b_i) x^i}.
\]

 
We may also define an action of $H$ on $PGL(2,F)$, as follows. Given $\sigma\in H$ and $s\in PGL(2,F)$, with $s$ expressed as a linear fractional transformation in the form $s(x)=(ax+b)/(cx+d)$, we define $\sigma.s$ by
\[
\sigma.s(x)=\frac{\sigma(a)x+\sigma(b)}{\sigma(c)x+\sigma(d)}.
\]
It is clear that $\sigma.s\in PGL(2,F)$ and furthermore, we have a homomorphism from $H$ into the automorphism group
of $PGL(2,F)$. This homomorphism is injective, as may be proved by elementary Galois theory.

Let $G$ be a finite subgroup of $PGL(2,F)$ and $H$ a subgroup of automorphisms of $F$, as above. Given $\sigma\in H$, we may define
another subgroup $\sigma.G$ by $\sigma.G = \{ \sigma.s | s\in G\}$. We say that $H$ fixes $G$ or that $G$ is $H$-invariant if
$\sigma.G=G$ for all $\sigma\in H$. We note that if $H$ fixes $G$, we inherit  a homomorphism from $H$ to the automorphism group
$Aut(G)$ of $G$. A little later in this section, we will give a sufficient condition for this homomorphism to be faithful. In these circumstances, the question of whether the automorphism of $G$ induced by a given element of $H$ is inner or outer is relevant, but it seems to be a delicate to give a definitive answer.

In connection with these questions, we need to make precise  a concept that we shall apply. As before, let $G$ be a finite subgroup of $PGL(2,F)$. We say that $G$ is defined over the subfield $K$ of $F$ if, for each element $s$ of $G$, where $s(x)=(ax+b)/(cx+d)$, there exists a nonzero $e$ in $\{a, b, c, d\}$ such that $a/e$, $b/e$, $c/e$ and $d/e$ are all in $K$. This means  that when the elements
of $G$ are expressed as fractional transformations, we can assume that all coefficients lie in $K$.
We also say that $G$ is not a subgroup of
$PGL(2,K)$ for any proper subfield $K$ of $F$ if $G$ is not defined over any such $K$. 

We may now proceed to prove a sufficient condition for a Galois group action to be faithful.

\begin{lemma} \label{faithful_Galois_action}
Let $F$ be a field and let $G$ be a finite subgroup of $PGL(2,F)$. Let $H$ be a nontrivial finite group of automorphisms of $F$.
Assume that $G$ is not a subgroup of $PGL(2,K)$ for any proper subfield $K$ of $F$. 
If $H$ fixes $G$ then $H$ acts faithfully
as a group of automorphisms of $G$.
\end{lemma}

\begin{proof}
Suppose that $\rho\in H$ acts trivially on $G$. Let $s$ be any element of $G$, with $s(x)=(ax+b)/(cx+d)$, where $a, \ldots, d$ are in $F$.
Then, applying $\rho$ to $s$, we obtain
\[
\rho.s(x)=s(x)=\frac{\rho(a)x+\rho(b)}{\rho(c)x+\rho(d)}=\frac{ax+b}{cx+d}.
\]
It follows that there is some nonzero element $\lambda$ in $F$ with
\[
\rho(a)=\lambda a,\qquad \rho(b)=\lambda b,\qquad \rho(c)=\lambda c,\qquad \rho(d)=\lambda d.
\]

Since $s$ is invertible, at least two of its entries are nonzero. For definiteness, we shall assume that $a$ is nonzero, but the argument
can trivially be modified if $a$ is zero, in which case some other entry is nonzero. We have $\lambda=a^{-1}\rho(a)$ and hence
\[
\rho(b)=a^{-1}\rho(a)b, \qquad \rho(c)=a^{-1}\rho(a)c, \qquad \rho(d)=a^{-1}\rho(a)d.
\]
We deduce that $a^{-1}b$, $a^{-1}c$, $a^{-1}d$ are all in the fixed field $K$, say, of $\rho$. Thus we have
\[
b=a\alpha,\qquad c=a\beta,\qquad d=a\gamma,
\]
where $\alpha, \beta, \gamma$ are in $K$. 

Since this holds for $s$ in $G$, we see that $G$ is a subgroup of $PGL(2,K)$. Given our hypothesis on $G$, $K$ equals $F$ and hence
$\rho=1$. This implies that $H$ acts faithfully as automorphisms of $G$, as required.
\end{proof}

We continue with the hypothesis that $G$ is a finite subgroup of $PGL(2,F)$ and $H$ is a finite group of automorphisms of $F$.
Recall that $G$ acts on $F(x)$ via $(s.r) (x)=r(s(x))$ for 
 $r\in F(x)$, $s\in G$. If $\sigma \in H$, then observe that
\begin{equation}\label{actions}
\sigma.(r(s(x)))=r(\sigma.s(x)).
\end{equation}

We want to relate the property that $H$ fixes $G$ to a property of the generator of the fixed subfield of $G$ acting on $F(x)$, as the next lemma informs us.

 \begin{lemma}\label{coeffs_in}
 Let $F$ be a field, and let $F(x)$ be the field of rational functions in $x$ over $F$.
 Let $G$ be any finite subgroup of $PGL(2,F)$.
 Assume that $G$ is not a subgroup of $PGL(2,K)$ for any proper subfield $K$ of $F$.
 Let $f(x)/g(x)$ be a nonconstant coefficient of the orbit polynomial $P_G(T)$ of $G$,
  with one of $f$, $g$ monic.
 Let $H$ be a nontrivial  finite group of automorphisms of $F$, with fixed field $E$.
 Then $H$ fixes $G$  if and only if  $f$ and $g$ have coefficients in $E$.
 \end{lemma}

\begin{proof}
Let $P_G(T)$ be the orbit polynomial of $G$.
The roots of $P_G(T)$ are the elements $s(x)$ for $s\in G$.
We consider the  $s(x)$ as elements of $F(x)$, as usual.

First  assume that $\sigma.G=G$ for all $\sigma \in H$.
Therefore the map $s(x)\mapsto (\sigma.s)(x)$ is a permutation of the roots of the orbit polynomial of $G$.
Any permutation of the roots of a polynomial fixes all the coefficients of the polynomial,
when the action is extended to the coefficients in the natural way.
We conclude that $H$ fixes the coefficients of the orbit polynomial.

Let $\Phi(x)=f(x)/g(x)$ be any nonconstant coefficient of $P_G(T)$.
We assume $g$ is monic (the same argument works if $f$ is monic).
Let $\sigma\in H$.
By the previous paragraph $\sigma.\Phi=\Phi$, so $\sigma.f/\sigma.g=f/g$.
This implies $(\sigma.f) g=(\sigma.g) f$, so $g$ divides $\sigma.g$ since $\gcd (f,g)=1$.
Since $g$ and $\sigma.g$ have the same degree, and $g$ is monic, we have $g=\sigma.g$,
and then also $f=\sigma.f$.
This is true for any $\sigma\in H$, and because the fixed field of $H$ is $E$,
the polynomials $f$ and $g$ must have coefficients in $E$.

For the converse, assume that $f$ and $g$ have coefficients in $E$.
Then $\Phi(x)=f(x)/g(x)$ has coefficients in $E$, so $\Phi$ is fixed by $H$.
Let $\sigma\in H$, and let $s\in G$.
By \eqref{actions} we have
$\sigma.(\Phi (s(x)))=\Phi(\sigma.s(x))$,
and since $\Phi$ is both $G$-invariant and fixed by $H$, this becomes
$\Phi(x)=\Phi(\sigma.s(x))$.
This means that $\Phi$ is invariant under the group $\sigma.G$ as well as the group $G$.
If $\sigma.G\not=G$ then $\Phi$ is invariant under the group generated by $\sigma.G$ and $G$,
which is larger than $G$, and this is a contradiction to the Galois correspondence.
Hence $\sigma.G=G$ for all $\sigma\in H$.
\end{proof}

We remark that in the examples considered in this section, we usually calculate the polynomials $f(x)$ and $g(x)$ using a computer algebra
system and then we can check if these polynomials lie in $E[x]$. Then we may deduce that $H$ fixes $G$. However, in a small example it is sometimes 
clear by inspection that $H$ fixes $G$, and then there is no 
obvious advantage in appealing to Lemma \ref{coeffs_in}.

In light of Lemma \ref{coeffs_in}, when $H$ fixes $G$ we may then ask for the Galois group of 
$f(T)-xg(T)$ over $E(x)$.
This question can be delicate, as the  examples below will  show. It is partly related to whether $H$ induces any inner automorphisms.

Firstly we state one obvious connection between the Galois groups over $F(x)$ and $E(x)$.
The following lemma is well known from elementary Galois theory.

\begin{lemma}\label{comp}
Let $L$ and $F$ be Galois extensions of $E$, contained in a fixed algebraic closure of $E$.
Then $Gal(LF/F) \cong Gal(L/L\cap F)$.
\end{lemma}

Next we note the following.

\begin{lemma}\label{spl_contain}
Let $F$ be a field, and let $F(x)$ be the field of rational functions in $x$ over $F$.
 Let $G$ be any finite subgroup of $PGL(2,F)$.
 Assume that $G$ is not a subgroup of $PGL(2,K)$ for any proper subfield $K$ of $F$.
  Let $f(x)/g(x)$ be a nonconstant coefficient of the orbit polynomial $P_G(T)$ of $G$.
 Let $L_F$ be the splitting field of $f(T)-xg(T)$ considered as an element of $F(x)[T]$.
 Let $H$ be a nontrivial finite group of automorphisms of $F$, with fixed field $E$, and suppose that $H$ fixes $G$.
  Assume that one of $f$, $g$ is monic. Then $f$ and $g$ have coefficients in $E$. Furthermore, 
 let $L_E$ be the splitting field of $f(T)-xg(T)$ considered as an element of $E(x)[T]$. 
 Then  $F(x)\subseteq L_E$
 and $L_E=L_F$.
\end{lemma}

\begin{proof}
We first observe that as $H$ fixes $G$, both $f$ and $g$ have coefficients in $E$. This follows from Lemma \ref{coeffs_in}.
Suppose that $F(x)$ is not contained in $L_E$.
Consider the following tower of fields.

  \setlength{\unitlength}{1mm}
   \begin{tabular}{c}

      \begin{picture}(62,67)(-31,0)
         \put(5,63){\makebox(0,0){$F(x) L_E=L_F$}}

         \put(0,60){\line(-2,-3){8}}
         \put(-12,44){\makebox(0,0){$F(x)$}}
         \put(-8,40){\line(2,-3){8}}

 \put(-6,33){\makebox(0,0){$G_1$}}
  \put(7,33){\makebox(0,0){$G_2$}}

         \put(0,60){\line(2,-3){8}}
         \put(10,44){\makebox(0,0){$L_E$}}
         \put(8,40){\line(-2,-3){8}}

         \put(0,23){\makebox(0,0){$F(x)\cap L_E$}}
          \put(0,19){\line(0,-3){11}}
          \put(0,3){\makebox(0,0){$E(x)$}}
      \end{picture}
 \end{tabular}
 
Let $G_1$ be the Galois group of $F(x)$ over $F(x)\cap L_E$, and let
$G_2$ be the Galois group of $L_E$ over $F(x)\cap L_E$.
Then the Galois group of the compositum $F(x) L_E=L_F$ over $F(x)\cap L_E$
is the direct product $G_1 \times G_2$.

Note that $G=Gal(L_F/F(x))$ by Theorem \ref{main1},
and  $G_2\cong G$ by Lemma \ref{comp}.
Also, $G_1$ is a  subgroup of $H$, since $E(x) \subseteq F(x)\cap L_E$
and $Gal(F(x)/E(x))\cong Gal(F/E)=H$.
And $G_1$ is a normal subgroup of $H$ because $F(x)$ and $L_E$ are both Galois
extensions of $E(x)$, so $F(x)\cap L_E$ is a Galois extension of $E(x)$.

As $G$ is not a subgroup of $PGL(2,K)$ for any proper subfield $K$ of $F$, $H$ acts faithfully on $G$, by Lemma \ref{faithful_Galois_action}.
It follows that we can choose $\tau \in G_1$
and $s\in G$ such that $\tau$ does not fix $s$.
Suppose $s$ is given by $s(x)=(ax+b)/(cx+d)$.

Let $\alpha$ be a  root of $f(T)-xg(T)$ in $L_F$.
Since $G$ acts transitively and regularly on the roots, 
$s(\alpha)= (a\alpha +b)/(c \alpha +d)$ is also a root and $s$ is the
unique element of $G$ that maps $\alpha$ to $s(\alpha)$.

Apply the automorphism $(\tau,1)$
to $s(\alpha)$, viewing $Gal(L_F/F(x)\cap L_E)$ as $G_1\times G_2$.
This element $(\tau,1)$ will fix $\alpha$ and $s(\alpha)$
and so we obtain  $s(\alpha) =  (\tau(a) \alpha +\tau(b))/(\tau(c) \alpha +\tau(d))$.
By uniqueness of $s$ we get
$(a\alpha +b)/(c \alpha +d)=(\tau(a) \alpha +\tau(b))/(\tau(c) \alpha +\tau(d))$
so there exists nonzero $\lambda$ such that
$\tau(a)=\lambda a$, $\tau(b)=\lambda b$, 
$\tau(c)=\lambda c$, and $\tau(d)=\lambda d$.
We now apply the proof of Lemma \ref{faithful_Galois_action}. At least one of $a$ and $b$ is nonzero, and we assume for definiteness
that $a$ is nonzero.  These equations imply 
that $\tau$ fixes $b/a$, $c/a$ and $d/a$.
Therefore $\tau$ fixes $s(x)=(x+b/a)/(cx/a+d/a)$,
and this is a contradiction to the choice of $\tau$ and $s$.
So $F(x)\subseteq L_E$.

Since $L_F$ is the compositum of $F(x)$ and $L_E$, and since $F(x)\subseteq L_E$, we
obtain  $L_E=L_F$.
\end{proof}

Next we demonstrate that in these circumstances, the Galois group over $E(x)$ 
is in fact a split extension of the Galois group over $F(x)$

\begin{lemma}\label{split}
Assume the hypotheses of Lemma \ref{spl_contain}.
Consider  $f(T)-xg(T)$ as an element of $E(x)$.
Let $L=L_F=L_E$.
Let $\Gamma = Gal(L/E(x))$.
Then  $\Gamma$ is a split extension of $G$.
\end{lemma}

\begin{proof}
By the Galois correspondence
$\Gamma / G \cong H$.
We will show that this extension of $G$ is always split.
Let $\alpha$ be a root of $f(T)-xg(T)$ in $L$.
Let $\Gamma_\alpha$ be the stabilizer of $\alpha$.
Then $G\cap \Gamma_\alpha = \{1\}$ because $G$ acts transitively and
regularly on the roots of $f(T)-xg(T)$.
Also, $|\Gamma_\alpha|=|\Gamma |/|G|$ by the orbit-stabilizer theorem.
Therefore
\[
|G \Gamma_\alpha | = |G| |\Gamma_\alpha| /|G\cap \Gamma_\alpha|=
 |G| |\Gamma_\alpha| = |\Gamma|
\]
and so $\Gamma = G \Gamma_\alpha$.
Thus $\Gamma$ is a semidirect product of $G$ and $ \Gamma_\alpha$.
\end{proof}

We note that
\[
H\cong \Gamma/G=G \Gamma_\alpha/G\cong
\Gamma_\alpha/(G\cap \Gamma_\alpha)\cong \Gamma_\alpha.
\]

We return to consider some examples in the light of the theory we have developed.

{\bf Example 1 continued:} 
Let $F=\mathbb{Q}(i)$ where $i^2=-1$.
 Then $PGL(2,F)$ contains a subgroup isomorphic to $A_4$ as outlined in Example 1 earlier.
By Theorem \ref{main1} 
\[
f_1(T)=T^{12} - 33T^8 - 33T^4+1-x(T^{10} - 2T^6 + T^2)
\]
is irreducible over $F(x)$ and has Galois group $A_4$.
By inspection or by Lemma \ref{coeffs_in}
we can see that $Gal(\Q(i)/\Q)$ permutes the elements of $A_4$.
We may ask for the Galois group of $f_1$ over $\Q(x)$.
By Lemma \ref{split},  $Gal(f_1/\Q (x))$ will be a split extension of $A_4$.
It turns out (using Magma)
that the Galois group of $f_1$ over $\Q(x)$ is the direct product $A_4 \times C_2$.



{\bf Example 2 continued:} 
Let $F=\mathbb{Q}(\omega)$ where $\omega^2+\omega+1=0$.
 Then $PGL(2,F)$ contains a subgroup isomorphic to $A_4$ as outlined in Example 2 earlier.
By Theorem \ref{main1} 
\[
f_2(T)=T^{12} + 264T^6 + 440T^3 + 24-x(T^9-3T^6+3T^3-1)
\]
is irreducible over $F(x)$ and has Galois group $A_4$.
By inspection or by Lemma \ref{coeffs_in}
we can see that $Gal(\Q(\omega )/\Q)$ permutes the elements of $A_4$.
By Lemma \ref{split},  $Gal(f_2/\Q (x))$ will be a split extension of $A_4$.
It turns out (using Magma) that the Galois group  of $f_2$ over $\Q(x)$ is $S_4$.


\bigskip

The appearance of two different split extensions of order 24 in these two examples, when the
original group is $A_4$ is both cases, demonstrates
that determining the Galois group exactly when descending to the subfield $E(x)$ is not trivial.

{\bf Example 3 continued:}
Let $F$ be any field of characteristic different from 2 in which 2 is a square, and let $c=1+\sqrt{2}$.
As shown earlier, 
\[
f_3(T)=T^{8} -28T^6 + 70T^4 - 28T^2 + 1 -x(T^7 - 7T^5 + 7T^3 - T)
\]
is irreducible over $F(x)$ and has Galois group $C_8$, a cyclic group of order 8.
Let us assume $F=\Q(\sqrt{2})$.
Since $f_3$ has coefficients in $\Q$,
by Lemma \ref{coeffs_in} we get that $Gal(\Q(\sqrt{2})/\Q)$ permutes the elements of $C_8$.
By Lemma \ref{split},  $Gal(f_3/\Q (x))$ will be a split extension of $C_8$.
The Galois group of $f_3$ over $\Q(x)$ has order 16, and according to Magma it
is the metacyclic or modular group of order 16 (not dihedral, semidihedral or dicyclic).

{\bf Example 5:}
From Klein's copy of $S_4$ (see \cite{K}) defined over $\Q(i)$ we obtain 
the polynomial
$T^{24} + 759T^{16} + 2576T^{12} + 759T^8 + 1 -xT^4(T^4-1)^4$
which has Galois group $S_4$ over $\Q(i)(x)$.
Over $\Q(x)$ 
it has Galois group $S_4 \times C_2$, according to Magma.

{\bf Example 6:}
From Klein's copy of $A_5$ (see \cite{K}) defined over $\Q(\epsilon)$ where $\epsilon$ is
a primitive fifth root of unity, for $f(T)-xg(T)$ over $\Q(x)$
we get a Galois group of order 240.
It is a semidirect product of $A_5$ and $C_4$, according to Magma.

\section{Galois groups of affine type}

Let $F$ be a field of prime characteristic $p$ and let $K$ be a Galois extension of $F$ of finite degree with
 Galois group $H$. We are interested in finding finite $H$-invariant $\mathbb{F}_p$-subspaces in $K$ (we may use $\mathbb{F}_q$-subspaces
 if we know that $\mathbb{F}_q$ is contained in $F$). We call such subspaces $\mathbb{F}_pH$-submodules.
 
 A simple way to proceed is as follows. Let $\alpha$ be an element in $K\setminus F$ and let $\Omega$ be the orbit of $\alpha$ under
 the action of $H$. Then the $\mathbb{F}_p$-span of the elements of $\Omega$ is an $H$-invariant subspace of dimension at most $|H|$. 
We call a subspace of this type a cyclic $\mathbb{F}_pH$-submodule and $\alpha$ a generator of the submodule.

A better way to show the existence of such $\mathbb{F}_pH$-submodules is to appeal to the normal basis theorem. Recall that if
$H=\{h_1, \ldots, h_m\}$, where $m=|H|$, there is an $F$-basis $e_{h_1}$, \dots, $e_{h_m}$ of $K$ over $F$ such that
\[
h(e_{h_i})=e_{hh_i}
\]
for all $h$ in $H$ and $1\leq i\leq m$. If we take the $\mathbb{F}_p$-span of this basis, we obtain an $H$-invariant subspace of dimension
$|H|$. This subspace is isomorphic to the regular $\mathbb{F}_pH$-module. It is a cyclic $H$-module, being generated by any element
$e_h$. Furthermore, any cyclic $\mathbb{F}_pH$-module is a homomorphic image of this regular module and conversely, 
any homomorphic image of it is cyclic. This enables us to construct (in principle, at least) a large number of cyclic $\mathbb{F}_pH$-submodules using $K$.

Now let $U$ be any finite dimensional $\mathbb{F}_pH$-submodule in $K$. For each element $u$ of $U$, we define a $K$-linear automorphism
$s_u$ of the function field $K(x)$ by setting
\[
s_u(x)=x+u.
\]
Since we have $s_us_w=s_{u+w}$ for any $w\in U$, it is clear that these elements form an elementary abelian $p$-subgroup of $K$-automorphisms of
$K(x)$ of order $|U|$, which we shall denote by $G_U$.

The Galois group $H$ acts on $K(x)$ by field automorphisms that
are $F$-linear, but not $K$-linear. The group $G_U$ is normalized by $H$ in its action on $K(x)$, since $U$ is $H$-invariant. Consider the element
\[
\Phi(x)=\prod_{u\in U}(x-u)
\]
in $K(x)$. It is $G_U$-invariant and, since it has degree $|U|=|G_U|$ when we consider it as a polynomial in $x$, it generates the subfield $K(x)^{G_U}$ of fixed points of $G_U$.

We now replace the variable $x$ by the more familiar variable $T$ and form the polynomial
\[
\Psi(T)=\Phi(T)-y,
\]
in $K(y)[T]$, $y$ being transcendental over $K$. The polynomial $\Psi$ has Galois group isomorphic to $G_U$ over $K(y)$. 
This follows from the earlier part of the paper. However, the coefficients
of $\Psi$ lie in $F(y)$, as $\Psi$ is invariant under the action of $H$, and
when we consider $\Psi$ as an irreducible polynomial in $F(y)[T]$, its splitting field has degree $|G_U||H|$ over $F(y)$. Its Galois group,
$\Gamma$, say, contains a normal subgroup, $\Gamma_1$, say, isomorphic to $G_U$. The quotient $\Gamma/\Gamma_1$ is isomorphic to $H$.
Furthermore, $\Gamma$ is a split extension of $\Gamma_1$ by a subgroup isomorphic to $H$. This follows from our work.

To elucidate the structure of $\Gamma$, we need to understand how a complement of $\Gamma_1$ acts by conjugation on $\Gamma_1$. This can be explained in terms of the roots of $\Psi$. Let $\alpha$ be a root of $\Psi$ in a splitting field over $K$. Then it is easy
to see that $K(y)(\alpha)$ is a splitting field for $\Psi$ over $K(y)$ and the roots of $\Psi$ are the elements $\alpha+u$, where $u$ runs
over the elements in $U$. 

Let $\Gamma_\alpha$ be the stabilizer of  $\alpha$ in the permutation action of $\Gamma$ on the roots of $\Psi$. Since $\Gamma_\alpha$ fixes
$\alpha$ and permutes the roots, it is straightforward to see that $\Gamma_\alpha$ maps $U$ into itself. Furthermore, $\Gamma_\alpha$ act
as field automorphisms of $U$ in the same manner as $H$. 

The subgroup $\Gamma_1$, which is a complement to $\Gamma_\alpha$, acts regularly on the roots, and fixes each element of $U$. If we label
the elements of $U$ as $u_i$, $1\leq i\leq |U|$, with $u_1=0$, we can correspondingly label the elements of $\Gamma_1$ as $s_i$, $1\leq i\leq |U|$, where 
\[
s_i(\alpha+u)=\alpha+u+u_i
\]
for all $u\in U$. Note then that
\[
s_i(s_j(\alpha))=\alpha+u_i+u_j=s_j(s_i(\alpha)).
\]

The next result, although elementary, enables us to understand better the internal structure of $\Gamma$.

\begin{thm} \label{isomorphism_theorem}
 With the notation previously introduced, $\Gamma_1$
 is isomorphic to $U$ when considered as an $\mathbb{F}_p\Gamma_\alpha$-module under conjugation.
\end{thm}

\begin{proof}
We define a mapping $\theta:\Gamma_1\mapsto U$ by $\theta(s_i)=u_i$ for all $i$. Since $s_is_j$ maps $\alpha$ to $\alpha+u_i+u_j$, it follows that
\[
\theta(s_is_j)=u_i+u_j
\]
and we deduce that $\theta$ is an $\mathbb{F}_p$-isomorphism.

Let $t$ be an arbitrary element of $\Gamma_\alpha$. We have then
\[
ts_it^{-1}(\alpha)=\alpha+t(u_i).
\]
It follows that 
\[
\theta(ts_it^{-1})=t(\theta(s_i))
\]
for all $i$. This proves that $\theta$ is an $\mathbb{F}_p\Gamma_\alpha$-isomorphism, as required.
\end{proof}

Note that as $\Gamma_\alpha$ is isomorphic to $H$, all $\mathbb{F}_p\Gamma_\alpha$-modules are $\mathbb{F}_pH$-modules. 

\begin{cor} \label{regular_module}
Let $F$ be a field of prime characteristic $p$ and suppose that $F$ has a Galois extension of finite degree with Galois group $H$. Then there exists a Galois extension of the function field $F(y)$ with Galois group $\Gamma$ having the following properties. $\Gamma$ contains a normal elementary abelian $p$-subgroup $\Gamma_1$ of order $p^{|H|}$, complemented by a subgroup isomorphic to $H$. As an $\mathbb{F}_pH$-module, $\Gamma_1$ is isomorphic to the regular $\mathbb{F}_pH$-module.
\end{cor}

We can extend this result slightly when we take account of the fact that cyclic modules are homomorphic images of the regular module.

\begin{cor} \label{cyclic_modules}
Assume the hypotheses of Corollary \ref{regular_module}. Let $V$ be a cyclic $\mathbb{F}_pH$-module. Then there
exists a Galois extension of the function field $F(y)$ with Galois group $\Delta$ having the following properties. $\Delta$ contains a normal elementary abelian $p$-subgroup $\Delta_1$ of order $|V|$, complemented by a subgroup isomorphic to $H$. As an $\mathbb{F}_pH$-module, 
$\Delta_1$ is isomorphic to $V$.
\end{cor}

\begin{proof}
We use the group $\Gamma$ described in Corollary \ref{regular_module}. Since $V$ is a homomorphic image of the regular module $\mathbb{F}_pH$-module, there is a normal subgroup $\Gamma_2$, say, of $\Gamma$ contained in $\Gamma_1$ such that $\Gamma_1/\Gamma_2$ is isomorphic to
$V$ as an $\mathbb{F}_pH$-module. The group $\Gamma/\Gamma_2$ is itself a Galois group over $F(y)$ by the Galois correspondence, and it has the required property.
\end{proof}

We remark that any $\mathbb{F}_pH$-module that contains a unique maximal $\mathbb{F}_pH$-submodule is cyclic. Thus, for example,
an irreducible $\mathbb{F}_pH$-module is cyclic as a module.

\noindent {\bf Example 7:} The polynomial $T^{24}+T+x$ has Galois group over $\mathbb{F}_2(x)$ isomorphic to the Mathieu group $M_{24}$. 
The roots of this polynomial span a twelve-dimensional $\mathbb{F}_2M_{24}$ module. This module contains an irreducible
eleven-dimensional submodule, called the Todd module. The elements of this module are the roots of the linearized polynomial
\[
T^{2048}+x^{64}T^{512}+x^8T^{16}+x^{16}T^8+T.
\]
The polynomial $F(x,T)/T$ is the product of irreducible polynomials of degree 276 and 1771. We presume that the roots of the degree 276 polynomial
are the sums of  two different roots of the original polynomial of degree 24. (Calculations done in Magma.)

Thus, if we denote this polynomial of degree 2048 by $F(x,T)$, the polynomial $F(x,T)+y$ in $\mathbb{F}_2(x,y)$ has Galois group
of order $2^{11}|M_{24}|$. The Galois group contains a normal elementary abelian subgroup of order $2^{11}$ that is isomorphic
to the eleven-dimensional Todd module for $M_{24}$. 

\section{Galois groups derived from linearized polynomials}

\noindent Let $F$ be a field of prime characteristic $p$.
We call a polynomial $P(T)$ of the form
 \[
\sum_{i=0}^n a_i T^{p^i}
 \]
 a linearized polynomial or $p$-polynomial over $F$. We assume here that $a_n\neq 0$, so that $P$ has degree $p^n$.
 
Let us assume additionally that $a_0\neq 0$.  In this case, the derivative $P'(T)=a_0$ and thus $P$ is separable. It follows
straightforwardly that the roots of $P$ in any splitting field over $F$ are a vector space of dimension $n$
 over $\mathbb{F}_p$. Thus, it follows the Galois group of $P$ over $F$ is a subgroup of the general
 linear group $GL(n,p)$ over $\mathbb{F}_p$. 
 
 The main result of this section is a simple construction that uses linearized polynomials to realize Galois groups that act doubly
 transitive on the roots of the associated polynomials. 
 
 We require the following lemma in the course of our proof. We extract its proof as a separate statement, but remark that it must be a standard result in the theory of polynomials.
 
\begin{lemma} \label{remains_irreducible_in_extension_field}
Let $K$ be a field and let $f$ be an irreducible polynomial in $K[T]$. Let $L$ be an extension field of $K$ of finite degree, with $[L:K]=m$.
Suppose that $m$ is relatively prime to $\deg f$. Then $f$ is irreducible over $L$.
\end{lemma}

\begin{proof}
Let $f_1$ be an irreducible factor of $f$ in $L[T]$ and let $\alpha$ be a root of $f_1$ in an extension field of $L$.
Clearly, $\alpha$ is also a root of $f$ and hence $[K(\alpha):K]=\deg f$.  We now  obtain
\[
[L(\alpha):L][L:K]=m\deg f_1=[L(\alpha):K(\alpha)][K(\alpha):K]=\deg f [L(\alpha):K(\alpha)].
\]
Thus $\deg f$ divides $m\deg f_1$ and since $\deg f$ is relatively prime to $m$, we deduce that $\deg f$ divides $\deg f_1$.
It follows that $f_1$ is a scalar multiple of $f$, as $\deg f_1\leq \deg f$, and this proves the lemma.
\end{proof} 

\begin{thm} \label{doubly_transitive} 
Let $F$ be a field of prime characteristic $p$ and let $f(T)$, $g(T)$ be linearized $p$-polynomials in $F[T]$, with
$p^n=\deg f>\deg g\geq 1$. Suppose that $f(T)/T$ and $g(T)/T$ are relatively prime. Let $\lambda\neq 0$ be an element of $F$ such that $f(T)+\lambda$ and $g(T)$ are also relatively prime. Then the Galois group $\Gamma$ of $P_\lambda(x,T)=f(T)+\lambda -xg(T)$ over $F(x)$ acts doubly transitively on the roots of the polynomial.
\end{thm}

\begin{proof}
Let us first note that $P_\lambda$ is separable (it has no repeated roots). For, let the coefficient of $T$ in $f(T)$ and $g(T)$ be
$a_0$, $b_0$ respectively. Then $P_\lambda'=a_0-xb_0$. If this expression is 0, then $a_0=b_0=0$. But this implies that $T^{p-1}$ divides
both $f(T)/T$ and $g(T)/T$, contrary to the relatively prime hypothesis. Thus, $P_\lambda$ has no repeated roots,
and likewise, $(f(T)-xg(T))/T$ has no repeated roots.

We set $E=F(x)$ and let $L$ be a splitting field for $P_\lambda$ over $E$. Let $\alpha$ be a fixed root of $P_\lambda$ in $L$
and let $\Gamma_\alpha$ be the stabilizer of $\alpha$ in $\Gamma$. Let $\beta$ be any other root of $P_\lambda$ different from $\alpha$.
Then since $f$ and $g$ are linearized, we  find that $\alpha-\beta$ is a root $(f(T)-xg(T))/T$. The polynomial $(f(T)-xg(T))/T$ has degree
$p^n-1$ and it is irreducible in $E[T]$ by Proposition 7.5.5 (b) of Cox. Since $p^n-1$ is relatively prime to $[E(\alpha):E]=p^n$,
$(f(T)-xg(T))/T$ is irreducible in $E(\alpha)[T]$, by Lemma \ref{remains_irreducible_in_extension_field}. 

The roots of $(f(T)-xg(T))/T$ are $\alpha-\beta_i$, where $\beta_i$ runs through the roots of $P_\lambda$ different from $\alpha$. 
We intend to show that $\Gamma_\alpha$ acts transitively on these roots. This is equivalent to showing that $\Gamma_\alpha$ acts
transitively on the $\beta_i$, and thus that $\Gamma$ is doubly transitive.

Let $\Omega$ be the $\Gamma_\alpha$-orbit containing some fixed root $\beta$. Consider the polynomial
\[
Q(T)=\prod_{\beta_i\in \Omega}(T-(\alpha-\beta_i)).
\]
This polynomial is fixed under the action of $\Gamma_\alpha$ and hence its coefficients lie in the fixed field of
$\Gamma_\alpha$, which is $E(\alpha)$. 

We have already seen that each element $\alpha-\beta_i$ is a root of $(f(T)-xg(T))/T$, which we now know is irreducible of degree
$p^n-1$ over $E(\alpha)$. It follows from the irreducibility that $Q(T)$ above is a scalar multiple
of $(f(T)-xg(T))/T$ and thus $|\Omega|=p^n-1$. This completes the proof.
\end{proof}

While we have proved that the Galois group $\Gamma$ acts doubly transitively on the roots, it is doubly transitive of a special type, namely
affine type. This means that $\Gamma$ contains a normal elementary abelian $p$-subgroup that acts transitively and regularly on the roots.
We prove this assertion below.

\begin{thm} \label{regular_normal_subgroup}
Assume the hypotheses of Theorem \ref{doubly_transitive}. Then the Galois group $\Gamma$ contains a normal elementary abelian $p$-subgroup
that acts transitively and regularly on the roots. The stabilizer of any root is isomorphic to the Galois group of
$(f(T)-xg(T))/T$ over $F(x)$.
\end{thm}

\begin{proof}
Let $\alpha$ be a fixed root of $P_\lambda$. As have noted, as $\beta$ runs over the roots of $P_\lambda$, the elements $\alpha-\beta$ are the roots of $f(T)-xg(T)$ and they form an $n$-dimensional vector space, $V$, say, of dimension $n$ over $\mathbb{F}_p$. Thus each root of
$P_\lambda$ is uniquely expressible as $v+\alpha$ for some $v\in V$.

The splitting field $L$, say, of $P_\lambda$ over $E=F(x)$ contains the elements of $V$. The subfield, $M$, say, generated over $E$ by $V$ is normal over $E$, since the polynomial $f(T)-xg(T)$ is an element of $E[T]$. Let $\Gamma_1$ be the subgroup of $\Gamma$ that fixes $M$ elementwise. By standard Galois theory, $\Gamma_1$ is normal in $\Gamma$ and $\Gamma/\Gamma_1$ is isomorphic to the
Galois group of $(f(T)-xg(T))/T$ over $E$.

We claim that no nonidentity element of $\Gamma_1$ fixes a root of $P_\lambda$. For suppose that $\sigma\in\Gamma_1$ fixes a root, which by transitivity we may assume to be $\alpha$. Certainly, $\sigma$ fixes each element of $V$, by definition. Since the roots of $P_\lambda$ are
expressible as $v+\alpha$, where $v\in V$, it follows that $\sigma$ fixes all roots, and hence is the identity.

The rest of the proof is routine from the elementary theory of doubly transitive groups. We briefly sketch some details.
Since $\Gamma_1$ is a nontrivial normal subgroup of a doubly transitive group, it acts transitively on the roots. As we have seen
that no nonidentity element of $\Gamma_1$ fixes a root, the action of $\Gamma_1$ is transitive and regular. 

Let $\tau$ be any element of $\Gamma_1$ and suppose that $\tau(\alpha)=u+\alpha$, where $u\in V$. Then we see that
for all $v\in V$,
\[
\tau(v+\alpha)=u+v+\alpha.
\]
It follows that 
\[
\tau^p(v+\alpha)=pu+v+\alpha=v+\alpha.
\]
It follows that $\tau^p=1$ and thus $\Gamma_1$ has exponent $p$. Furthermore, if $\tau'\in\Gamma_1$ and $\tau'(\alpha)=w+\alpha$,
where $w\in V$, we have
\[
\tau'\tau(v+\alpha)=u+w+v+\alpha=\tau\tau'(v+\alpha).
\]
This implies that $\tau$ and $\tau'$ commute, and we deduce that $\Gamma_1$ is abelian. Hence $\Gamma_1$ is an elementary abelian
$p$-group, as claimed.
\end{proof}

The structure of $\Gamma$ depends to a large extent on the structure of the Galois
group of $(f(T)-xg(T))/T$, and this may depend on properties of the field $F$.

As far as choices of linearized polynomials $f$ and $g$ are concerned, there seems to be a lot of flexibility available. For example,
we may take $f(T)$ to be any linearized polynomial for which the coefficient of $T$ is nonzero, and then take $g(T)=\mu T^{p^m}$,
where $\mu$ is any nonzero element of $F$ and $m\geq 1$. 

Suppose also that we have linearized polynomials $f$ and $g$ such that $f(T)/T$ and $g(T)/T$ are relatively prime. Then it is easy to check
that if $\deg g=m$, there are at most $m$ elements $\lambda$ of $F$ such that $f(T)+\lambda$ and $g(T)$ have a nontrivial common factor.
Thus, for example, if $F$ is infinite, almost all elements of $F$ may be used.

\noindent{\bf Example 8:} Let $q$ be a power of the prime $p$ and $n$ be a positive integer. Consider the linearized polynomials 
\[
f(T)=T^{q^n}-T^q, \quad g(T)=T
\]
in $
\mathbb{F}_q[T]$. These satisfy the hypotheses of Theorem \ref{doubly_transitive}. Abhyankar has shown that the Galois group
of $(f(T)-xg(T)/T$ over $\mathbb{F}_q(x)$ is isomorphic to the general linear group $GL(n,q)$. If $\lambda$ is any non-zero element
of $\mathbb{F}_q$, the Galois group of $f(T)+\lambda-xg(T)$ over $\mathbb{F}_q(x)$ contains a normal elementary abelian subgroup
of order $q^n$ with quotient isomorphic to $GL(n,q)$. This group must be isomorphic to the affine general linear group.

We consider another rather similar construction that uses linearized polynomials to obtain Galois groups of affine type. This construction is fairly well known in the context of arithmetic and geometric monodromy groups of polynomials but it works in an especially satisfying way with
linearized polynomials. Let $F$ be a field of prime characteristic $p$ and let $\overline{F}$ denote the algebraic closure of $F$. Let
$x$ be transcendental over $\overline{F}$ and let $\overline{F}(x)$ denote the function field over $\overline{F}$. We will be interested
in subfields of $\overline{F}(x)$ of the form $E(x)$, where $E$ is an extension of $F$ of finite degree, considered as a subfield
of $\overline{F}$. 

Let $f(T)$ be a linearized polynomial in $F[T]$ of degree $p^n$. We make the assumption that the coefficient of $T$ in $f(T)$ is non-zero.
This ensures that $f(T)$ has no repeated roots. We consider the polynomial $f(T)+x$ in various subfields
of $\overline{F}(x)$.

\begin{lemma} \label{geometric_monodromy_group}
Let $f(T)$ be a linearized polynomial of degree $p^n$ in $F[T]$ and let $E$ be a splitting field for $f$ over $F$. Then the Galois group
of $f(T)+x$ over $E(x)$ is elementary abelian of order $p^n$. 
\end{lemma}

\begin{proof}
Let $\alpha$ be a fixed root of $f(T)+x$ in some splitting field over $E(x)$. Since $f(T)$ is linearized, it is clear
that the roots of $f(T)+x$ are the elements $\alpha+\lambda$, where $\lambda$ runs over the $p^n$ different roots of $f$ contained
in $E$. 

Let $G$ denote the Galois group of $f(T)+x$ over $E(x)$ and let $\sigma$ be an element of $G$. Then $\sigma(\alpha)=\alpha+\lambda$, where
$\lambda$ is a root of $f$. Let $\beta=\alpha+\mu$ be any other root of $f(T)+x$, where $\mu$ is some root of $f$. Since $\mu\in E$
and $\sigma$ fixes $E$ elementwise, we have
\[
\sigma(\beta)=\sigma(\alpha+\mu)=\alpha+\lambda+\mu=\beta+\lambda.
\]
We deduce in a straightforward way that
\[
\sigma^p(\beta)=\beta+p\lambda=\beta.
\]
This argument shows that $\sigma^p$ fixes all the roots of $f(T)+x$ and hence is the identity. We have thus shown that
$G$ has exponent $p$.

Let $\tau$ be any other element of $G$ and suppose that $\tau(\alpha)=\alpha+\kappa$, where $\kappa$ is some root of $f$. Then, as we have shown above, $\tau(\beta)=\beta+\kappa$ for any root $\beta$ of $f(T)+x$. It now follows that
\[
\tau\sigma(\beta)=\beta+\lambda+\kappa=\sigma\tau(\beta)
\]
and we deduce that $\tau\sigma=\sigma\tau$. Thus $G$ is an elementary abelian $p$-group. 

The polynomial $f(T)+x$ is irreducible of degree $p^n$ over $E(x)$ and
 thus its Galois group acts transitively on the roots. This implies
that $p^n$ divides $|G|$. However, it is clear that no nonidentity of $G$ fixes a root and thus $G$ acts regularly. As a consequence,
$|G|=p^n$, as required.
\end{proof}

It is clear that the same argument shows that the Galois group of $f(T)+x$ over $\overline{F}(x)$ is also
elementary abelian of order $p^n$. Thus we can say that the geometric monodromy group of $f(T)+x$ is the same elementary abelian group.

We proceed to determine the Galois group of $f(T)+x$ over the base field $F(x)$.

\begin{thm} \label{arithmetic_monodromy_group}
Assume the hypotheses of Lemma \ref{geometric_monodromy_group}. Let $\Gamma$ be the Galois group of $f(T)+x$ over $F(x)$. Then $\Gamma$ contains a normal subgroup $\Gamma_1$ isomorphic to the Galois group of $f(T)+x$ over $E(x)$. The quotient group
$\Gamma/\Gamma_1$ is isomorphic to the Galois group $G$ of $f(T)$ over $F$. If $\alpha$ is any root of $f(T)+x$ in its splitting field
over $F(x)$, and $\Gamma_\alpha$ is the stabilizer of $\alpha$, $\Gamma_\alpha$ is isomorphic to $G$ and $\Gamma$ is the semidirect
product of $\Gamma_1$ with $\Gamma_\alpha$.
\end{thm}

\begin{proof}
Let $L$ be a splitting field for $f(T)+x$ over $F(x)$ and let $\alpha$ be a root of the polynomial in $L$. We have shown
in the proof of Lemma \ref{geometric_monodromy_group} that the roots of $f(T)+x$ are the elements $\alpha+\lambda$, where
$\lambda$ runs over the roots of $f$. It follows that $L$ contains $E$. Since $L$ also contains $F(x)$ and $\alpha$, $L$ contains
$E(x)(\alpha)$, which is the splitting field for $f(T)+x$ over $E(x)$. We deduce from the minimality of splitting
fields that $L=E(x)(\alpha)$ and the splitting field for the polynomial over $F(x)$ is the same as the splitting field
over $E(x)$. 

Since $E(x)$ is clearly a normal extension of $F(x)$ contained in $L$, it is invariant under the action of $\Gamma$. Thus, restriction 
to $E(x)$ defines
a homomorphism of $\Gamma$ into the Galois group of $E(x)$ over $F(x)$. This homomorphism is surjective by basic Galois theory.
The kernel of the homomorphism is the group of $E(x)$-automorphisms of $L$, which Lemma \ref{geometric_monodromy_group} shows
is elementary abelian of order $p^n$.

Finally, we have shown that $\Gamma_1$ acts transitively and regularly on the roots and thus $\Gamma_1\cap \Gamma_\alpha=1$. Since $\Gamma_\alpha$ is the stabilizer
of a root in a transitive action, $|\Gamma:\Gamma_1|=p^n$. It follows that $\Gamma=\Gamma_1\Gamma_\alpha$, as required.
\end{proof}

Given this analysis of the Galois group $\Gamma$ of $f(T)+x$, we would like now to draw attention to aspects of the structure of
$\Gamma$ which encode information about how the polynomial $f(T)$ factors into irreducibles in $F[T]$. The observations we make are simple consequences of an identification of the subspace of roots of $f(T)$ and the normal subgroup $\Gamma_1$ of $\Gamma$. Of course, much of this
identification is already implicit in our proof of Theorem \ref{arithmetic_monodromy_group}. We start with further explanations.

The set of roots of $f(T)$ in the splitting field $E$ over $F$ is a vector space of dimension $n$ over $\mathbb{F}_p$, which we shall denote
by $V$ and call the subspace of roots. The Galois group $G$ of $f(T)$ acts faithfully on $V$ and thus $V$ is a faithful $\mathbb{F}_pG$-module.

In the group $\Gamma$, the normal elementary abelian $p$-subgroup $\Gamma_1$ may be identified with a vector space of dimension $n$ over
$\mathbb{F}_p$. The subgroup $\Gamma_\alpha$ of $\Gamma$ acts on $\Gamma_1$ by conjugation and we may thus consider $\Gamma_1$ to be an
$\mathbb{F}_p\Gamma_\alpha$-module. Furthermore, we know that $\Gamma_\alpha$ is isomorphic to $G$. Thus, it is not unreasonable to anticipate
that the following isomorphism theorem is true. 

\begin{thm} \label{isomorphism_theorem_for_modules}
With the notation previously introduced, the subspace $V$ of roots of $f(T)$, considered as an $\mathbb{F}_pG$ under Galois action, is isomorphic to $\Gamma_1$, when considered as as
an $\mathbb{F}_p\Gamma_\alpha$-module under conjugation. ($G$ and $\Gamma_\alpha$ are isomorphic.)
\end{thm}

\begin{proof}
Let $\alpha$ be a fixed root of $f(T)+x$ in a splitting field over $F(x)$. Given $\lambda$ in $V$, we have shown in Theorem
\ref{arithmetic_monodromy_group} that there is a corresponding element $s_\lambda$ in $\Gamma_1$ that is specified by its action on $\alpha$:
\[
s_\lambda(\alpha)=\alpha+\lambda.
\]
We define a mapping $\theta$ from $V$ to $\Gamma_1$ by setting
\[
\theta(\lambda)=s_\lambda.
\]
Since we have already shown that $s_\lambda s_\mu=s_{\lambda+\mu}$
for $\mu\in V$, $\theta$ is an isomorphism of abelian groups. 

Note that as $V$ is an $\mathbb{F}_p$-subspace of the splitting field of $f(T)+x$, $\Gamma_\alpha$ acts on $V$ in an obvious way. Let
$g$ be an arbitrary element of $\Gamma_1$. Then we have
\[
gs_\lambda g^{-1}(\alpha)=g(\alpha+\lambda)=\alpha+g(\lambda).
\]
This implies that $gs_\lambda g^{-1}=s_{g(\lambda)}$ and this gives us the required isomorphism with compatible group actions.
\end{proof}

the next result shows how properties of the conjugacy classes of $\Gamma$ contained in $\Gamma_1$ reflect the way in which
the polynomial $f(T)$ factors into irreducible divisors of in $F[T]$.

\begin{thm} \label{conjugacy_class_sizes}
There is a one-to-one correspondence between the conjugacy classes of the Galois group $\Gamma$ of $f(T)+x$ over $F(x)$ that are contained
in the normal subgroup $\Gamma_1$ and the monic irreducible factors of $f(T)$ that lie in $F[T]$.
Under the correspondence, the size of a conjugacy class is equal to the degree of the corresponding
irreducible factor. Thus, the number of conjugacy classes contained in $\Gamma_1$ equals the number of irreducible factors of $f(T)$
and this number is equal to the rank of $\Gamma$ considered as a permutation group on the roots of $f(T)+x$.

Furthermore, $\Gamma_1$ is a minimal normal subgroup of $\Gamma$ if and only if $f(T)$ has no nontrivial linearized polynomial divisor  other than $T$.
\end{thm}

\begin{proof}
Most of the proof consists of interpreting the consequences of Theorem \ref{isomorphism_theorem_for_modules} in two ways and  it is thus
 largely formal. Consider an irreducible factor of $f(T)$ in $F[T]$. There corresponds a $G$-orbit of roots of this polynomial
and consequently a $G$-orbit on the subspace of roots $V$. Theorem \ref{isomorphism_theorem_for_modules} shows that
this orbit on $V$ corresponds to an orbit of $\Gamma_\alpha$ acting on $\Gamma_1$ by conjugation and hence a conjugacy class of $\Gamma$
contained in $\Gamma_1$. The sizes of the two orbits are the same, and hence the degree of the irreducible factor of $f(T)$ is the size of the
conjugacy class.

The group $\Gamma$ acts transitively on the roots of $f(T)+x$ and the rank of $\Gamma$ is by definition the number of
orbits of the stabilizer $\Gamma_\alpha$ of a root $\alpha$ acting on the roots. By what we have just shown, 
the number of such orbits is the number irreducible factors of $f(T)$. 

Finally, suppose that $\Gamma_1$ is not minimal normal, but contains say a proper normal subgroup, $\Gamma_2$, say. Then $\Gamma_2$ is a union
of conjugacy classes of $\Gamma$ and is obviously $\Gamma_\alpha$ invariant. By Theorem \ref{isomorphism_theorem_for_modules},
there corresponds a subspace, $W$, say of the subspace $V$ of roots of $f(T)$ that is $G$-invariant.
Then the polynomial
\[
f_1(T)=\prod_{w\in W} (T-w)
\]
is linearized and it divides $f(T)$, whose roots are all the elements of $V$. Finally, $f_1(T)$ is fixed by the Galois group
$G$, since its roots are a $G$-orbit, and thus lies in $F[T]$. Conversely, a nontrivial linearized factor of $f(T)$ in $F(T)$ determines
a subspace of $V$ that is $G$-invariant, and our correspondence theorem shows that this is equivalent to a $\Gamma_\alpha$-invariant proper
subgroup of $\Gamma_1$, which is then normal in $\Gamma_1$.

\end{proof}

\smallskip

{\bf Example:} Suppose that we take $F=\mathbb{F}_p$, with $p$ a prime, and let $n$ be a positive integer. Let $f(T)=T^{p^n}-T$ in
$\mathbb{F}_p[T]$. Then the Galois group of $f(T)$ over $\mathbb{F}_p$ is cyclic of order $p$. The Galois group $\Gamma$ of $f(T)+x$
over $\mathbb{F}_p(x)$ has order $n\cdot p^n$. $\Gamma$ contains a normal elementary abelian $p$-subgroup $\Gamma_1$ of order $p^n$, which is complemented by a cyclic group of order $n$. We can identify $\Gamma$ with the regular wreath product of a cyclic group of order $p$
by a cyclic group of order $n$. In the wreath product, $n$ commuting copies of a group of order $p$ are transitively permuted by conjugation action.
The number of conjugacy classes of $\Gamma$ contained in $\Gamma_1$ is the number of irreducible monic polynomials in $\mathbb{F}_p[T]$ of degree dividing $n$.

\section{Identifying groups obtained by descent}\label{more}

We want in this section to elucidate the structure of the Galois group $\Gamma$ obtained by descent, as described in Lemma \ref{spl_contain}, and relate
this structure to the groups $G$ and $H$. This of course is dependent on how $H$ acts on $G$. We have a large set of hypotheses that we need to
assume, and thus to save ourselves from an excess of repetition, we make an initial statement of what we are assuming, together with some consequences of our hypotheses that we have already proved. Our subsequent results will then be stated in terms of these hypotheses.

\noindent{\bf Hypotheses (H).} $F$ is an arbitrary field and $G$ is a finite subgroup of $PGL(2,F)$. We assume that $G$ is not a subgroup of $PGL(2,K)$
for any proper subfield $K$ of $F$. We have a nonconstant coefficient $f(x)/g(x)$ of the orbit polynomial of $G$, with one of $f$ and $g$ monic.
We also have  a finite group $H$ of automorphisms of $F$, with fixed field $E$, that fixes $G$ in its action on $PGL(2,F)$. (Thus both $f$ and $g$ have coefficients in $E$, by Lemma \ref{coeffs_in}.)

The polynomial $P$ defined by $P=f(T)-xg(T)$ is an element of $E(x)[T]$ and $\Gamma$ is the Galois group of $P$ over $E(x)$. (Thus
$\Gamma$ contains a normal subgroup $\Gamma_1$, say, isomorphic to $G$, which fixes $F(x)$ elementwise.
$\Gamma_1$ is the Galois group of $P$ over $F(x)$.) $\Lambda$ is the set of roots
of $P$ in a splitting field and $\alpha$ is an element of $\Lambda$. $\Gamma_\alpha$ is the stabilizer of $\alpha$ in $\Gamma$. (Then
$\Gamma_\alpha$ is a complement for $\Gamma_1$ in $\Gamma$ and $\Gamma_\alpha$ is isomorphic to $H$, by Lemma \ref{split} and the subsequent
observations.)

Our main objective will be to find conditions to guarantee that the centralizer of $\Gamma_1$ in $\Gamma$ is contained in $\Gamma_1$ (and thus
is the centre of $\Gamma_1$). When we can be sure that the centralizer is in $\Gamma_1$, we may be able to identify $\Gamma$ in terms
of the automorphism group of $\Gamma_1$. We will illustrate this principle with some naturally occurring examples of groups $G$ 
and $H$. 

We continue with some further observations about our hypotheses. Let $L$ be a splitting field for $P$ over $E(x)$. We have proved that
$L$ contains $F$. Since $\Gamma$ fixes $E$ elementwise, and $F$ is a Galois extension of $E$, it follows that $\Gamma$ maps $F$ into itself
and we thus have a homomorphism from $\Gamma$ into the Galois group $H$ of $F$ over $E$. $\Gamma_1$ is in the kernel of the homomorphism, as it acts trivially on $F(x)$, and we have already shown that the homomorphism is surjective. It follows that in its action on
$F$, $\Gamma_\alpha$ induces the full Galois group $H$ and we will consequently identify $\Gamma_\alpha$ with $H$. 

We can now formulate a simple general principle concerning the number of roots fixed by elements of $\Gamma_\alpha$. 

\begin{lemma} \label{number_of_roots_fixed}
Assume the hypotheses (H). Let $\tau$ be an element of $\Gamma_\alpha$. Then the number of elements of $\Lambda$ fixed by
$\tau$ equals the order of the subgroup of $G$ fixed by the automorphism of $F$ defined by $\tau$.
\end{lemma}

\begin{proof}
Let $\beta$ be an element of $\Lambda$ fixed by $\tau$. The proof of Theorem \ref{main1} shows that each element of $\Lambda$ has the form $\beta =s(\alpha)$ for a unique element $s$ of $G$. Let $s$ be defined as the fractional transformation $x\mapsto (ax+b)/(cx+d)$, where the coefficients $a$, \dots, $d$ are in $F$. We thus have
\[
\beta=\tau(\beta)=(\tau(a)\alpha+\tau(b))/(\tau(c)\alpha+\tau(d))=(a\alpha+b)/(c\alpha+d).
\]
Arguing as in the proof of  Lemma \ref{faithful_Galois_action},  we deduce  that $s$ is fixed by $\tau$ in its action on $G$, and conversely, any element of $G$ fixed by
$\tau$ determines a root fixed by $\tau$. Thus the number of roots fixed by $\tau$ is the order of this fixed subgroup of $G$.
\end{proof}

We proceed to the proof of an important technical lemma.

\begin{lemma} \label{pinning_down_the_centralizer}
Assume the hypotheses (H) and let $C$ denote the centralizer of $\Gamma_1$ in $\Gamma$. Suppose that $C$ is not contained in $\Gamma_1$. Then there is a nonidentity element $\tau$ of $\Gamma_\alpha$ such that the number of elements of $\Lambda$ fixed by $\tau$ is the order of the centralizer in $\Gamma_1$ of some nonidentity element of $\Gamma_1$.
\end{lemma}

\begin{proof}
We first note that as $\Gamma_1$ is normal in $\Gamma$, $C$ is also normal in $\Gamma$. Let $z$ be an element of $C$ not contained in $\Gamma_1$. Since $\Gamma_\alpha$ is a complement of $\Gamma_1$ in $\Gamma$, we can write
\[
z=g^{-1}\tau,
\]
where $g$ is an element of $\Gamma_1$ and $\tau$ an element of $\Gamma_\alpha$. We note that $\tau\neq 1$, since otherwise $z$ is an element
of $\Gamma_1$, contrary to hypothesis.

Let us now show that $g\neq 1$. For if $g=1$, it follows that $z$ is a nonidentity element of $\Gamma_\alpha$ and hence fixes $\alpha$. We derive a contradiction as follows. Let $h$ be any element of $\Gamma_1$. Then we have
\[
zh(\alpha)=hz(\alpha)=h(\alpha),
\]
since $z$ and $h$ commute. But $\Gamma_1$ acts transitively on $\Lambda$ and hence $z$ fixes all roots. This is the desired contradiction,
since only the identity of $\Gamma_1$ fixes all roots. Thus we have shown that $g\neq 1$.

Clearly, $\tau=gz=zg$ and we will proceed to calculate how many roots are fixed by $\tau$ when we represent $\tau$ in this form.
We claim that the answer is $|D|$, where $D$ is the centralizer of $g$ in $\Gamma_1$. Let us note that any element $d$ of $D$ commutes
with $zg=\tau$, since $z$ centralizes $\Gamma_1$. Thus the argument above shows that $\tau$ fixes $d(\alpha)$. Since $\Gamma_1$ acts regularly
on $\Lambda$, it follows that we have found $|D|$ different roots fixed by $\tau$. 

To complete the enumeration, we show that any root $\beta$ fixed by $\tau$ has the form $e(\alpha)$ for some $e$ in $D$. To achieve this, 
we recall that $\beta=w(\alpha)$ for a unique $w$ in $\Gamma_1$. Thus we obtain
\[
zgw(\alpha)=w(\alpha).
\]
Since $w$ commutes with $z$, 
\[
zw^{-1}gw(\alpha)=\alpha=zg(\alpha).
\]
We cancel $z$ from the equality above and deduce that $g^{-1}w^{-1}gw$ fixes $\alpha$. Hence $w^{-1}gw=g$, since $\Gamma_1$ acts regularly, and thus $w$ centralizes $g$. This establishes the desired formula for the number of fixed roots.
\end{proof}

We may rapidly deduce what we consider to be a significant consequence of Lemma \ref{pinning_down_the_centralizer}.

\begin{thm} \label{abelian_subgroup}
Assume the hypotheses (H), with the further hypothesis that $G$ is abelian. Then $\Gamma_1$ is its own centralizer in $\Gamma$ and hence
$\Gamma/\Gamma_1$ acts faithfully on $\Gamma_1$ as a group of automorphisms.
\end{thm}

\begin{proof}
Let $C$ denote the centralizer of $\Gamma_1$ in $\Gamma$. Since $G$ is abelian, $\Gamma_1$ is contained in $C$. Suppose if possible that
$C$ is strictly larger than $\Gamma_1$. Then Lemma \ref{pinning_down_the_centralizer} states that there is a nonidentity
element $\tau$ of $\Gamma_\alpha$ such that the number of roots fixed by $\tau$ is the order of the centralizer in $\Gamma_1$ of an
element of $\Gamma_1$. However, as $\Gamma_1$ is abelian, this number is $|\Gamma_1|$. Thus $\tau$ fixes all roots, which we know is impossible. Therefore, $C=\Gamma_1$, as required.
\end{proof}



\section{An application in Characteristic 0}\label{dihedral_sect}

When $F$ is a field of characteristic zero, apart from the three groups related to the regular polyhedra, there are only cyclic and dihedral
subgroups of $PGL(2,F)$. We will apply the theory developed so far to investigate a Galois group related to the dihedral group. We 
take $F=\Q(\zeta_n)$ to be a cyclotomic field, generated by the primitive $n$-th root of unity $\zeta_n$. In this case, $PGL(2,F)$ contains a dihedral subgroup, $G$, say, of order $2n$ generated by the fractional
transformations
\[
s(x)=\zeta_n x,\quad t(x)=1/x,
\]
where $s$ has order $n$, $t$ has order 2, and $t$ inverts $s$. The polynomial $x^n$ is invariant under $s$ and we easily see that
\[
x^n+x^{-n}=(x^{2n}+1)/x^n
\]
generates the subfield of $G$-invariant rational functions in $F(x)$. 

It follows from Theorem \ref{main1} that the polynomial
\[
D(T)=T^{2n}-xT^n+1
\]
has Galois group $G=D_n$ over $F(x)$. We want to describe the Galois group of $D(T)$ over $\mathbb{Q}(x)$. It has order $2n\phi(n)$, where
$\phi$ denotes the Euler function, and we may anticipate that the group is related to the automorphism group of the dihedral group.
We intend to show that, when $n$ is odd, the Galois group is isomorphic to the direct product of the automorphism group of $D_n$ with a group of order 2. We omit a description of the Galois group when $n$ is even, not least because the answer is somewhat messy.

We start our analysis by determining the automorphism group of the dihedral group when $n$ is odd.

\begin{lemma} \label{automorphism_group_of_dihedral}
Let $G$ be a dihedral group of order $2n$, where $n$ is odd. Then the automorphism group of $G$ has order $n\phi(n)$ and 
it contains a normal cyclic subgroup of order $n$ generated by an inner automorphism.
\end{lemma}

\begin{proof}
We take $G$ to be generated by elements $s$ and $t$, as above, where $s$ has order $n$ and $t$ is an involution that inverts $s$.
Let $H$ be the cyclic subgroup of order $n$ generated by $s$. It is straightforward to see that $H$ is characteristic in $G$ and hence
we have a homomorphism from the automorphism group of $G$ into the automorphism group of $H$, obtained by restriction
to $H$. We shall now determine
the kernel $K$, say, of this homomorphism.

Let $\sigma$ be an element of $K$. Then we have $\sigma(s)=s$. Since $\sigma(t)$ is an involution
in $G$, it can be written as $s^at$ for some integer $a$. We claim that $\sigma$ is inner. For, since $|H|$ is odd, we can write
$s^a=s^{2b}$ for some integer $b$. Then we have
\[
s^{b}ts^{-b}=s^{2b}t=s^a t=\sigma(t).
\]
Since conjugation by $s^b$ fixes $H$ elementwise, $\sigma$ is the inner automorphism determined by $s^b$, as claimed.

Conversely, any inner automorphism of $G$ defined by an element of $H$ fixes $H$ elementwise and hence is in $K$. This implies that
$K$ is cyclic of order $n$ and is identical with the subgroup of inner automorphisms determined by $H$. We proceed to show that
any automorphism of $H$ may be extended to $G$, so that the homomorphism of one group to the other is surjective. Let
$m$ be an integer relatively prime to $n$ and define $\tau:G\to G$ by
\[
\tau(s^r)=s^{rm}, \quad \tau(s^rt)=s^{rm}t,
\]
for $1\leq r\leq n$. It is easy to verify that $\tau$ is an automorphism of $G$ that extends the automorphism $s\mapsto s^m$ of $H$.
Since there are $\phi(n)$ such automorphisms, we have proved the required result.
\end{proof}

\begin{thm} \label{dihedral_automorphism_group}

Let $n>1$ be an odd integer. Then the Galois group of $D(T)=T^{2n}-xT^n+1$ over $\mathbb{Q}(x)$ is isomorphic to the direct product
of the automorphism group of a dihedral group of order $2n$ with a group of order $2$.
\end{thm}

\begin{proof}
We again adopt  the notation of Hypotheses (H), with $F=\mathbb{Q}(\zeta_n)$, $E=\mathbb{Q}$, $G$ a dihedral group of order $2n$ and $H$ abelian of order $\phi(n)$.
  $\Gamma$ contains a normal subgroup $\Gamma_1$  isomorphic to the dihedral group, with quotient of order $\phi(n)$.
  
  Let $C$ denote the centralizer of $\Gamma_1$ in $\Gamma$. We certainly have $\Gamma_1\cap C=1$, since the centre of a dihedral group
  is trivial when $n$ is odd. On the other hand, $C$ is nontrivial. For if $C=1$, $\Gamma$ embeds into the automorphism group
  of $\Gamma_1$, which we know from Lemma \ref{automorphism_group_of_dihedral} has order $n\phi(n)$, whereas $|\Gamma|=2n\phi(n)$.
  This contradiction implies that $C>1$. We need to show that $|C|=2$.
  
  Let $z$ be a nonidentity element of $C$. As we observed in Lemma \ref{pinning_down_the_centralizer}, we can write $zg=\tau$, where
  $g$ and $\tau$ are nonidentity elements of $\Gamma_1$ and $\Gamma_\alpha$, respectively. The number of roots of $D(T)$ fixed by
  $\tau$ is the order of the centralizer of $g$ in $\Gamma_1$. Now as $\Gamma_1$ is isomorphic to a dihedral group of order $2n$, where $n$ is odd, there are two possibilities. Either $g$ has odd order, dividing $n$, in which case its centralizer in $\Gamma_1$ has order $n$, or $g$ has order 2, in which case its centralizer in $\Gamma_1$ has order 2. We claim that the first case cannot occur. For if the centralizer
  of $g$ has order $n$, $\tau$ fixes exactly $n$ roots of $D(T)$. However, $D(T)$ clearly has the property that if $\beta$ is a root, so also is its reciprocal. Since $\beta$ and $\beta^{-1}$ are different, $\tau$ must fix an even number of roots. This verifies our claim that
  $g$ has order 2.
  
  Since $z$ and $g$ commute, we may square the equation $zg=\tau$ to obtain $z^2=\tau^2$, and deduce that $z^2$ is in $\Gamma_\alpha$ and hence fixes the root $\alpha$. If $z^2$ is not the identity, the argument used in the proof of Lemma \ref{pinning_down_the_centralizer} implies that $z^2$ fixes all roots, since it centralizes $\Gamma_1$. This cannot occur and hence $z^2=1$. Thus we have shown
  that all nonidentity elements of $C$ have order 2. 
  
  Suppose if possible that $C$ has order greater than 2, and let $z_1$ be another element of order 2 in $C$ different from $z$. 
  Then we can write $z_1g_1=\tau_1$, where $g_1$ is an element of order 2 in $\Gamma_1$, and $\tau_1$ an element of order 2 in $\Gamma_\alpha$. When we multiply this equation by $zg=\tau$, and recall that $C$ centralizes $\Gamma_1$, we obtain
  \[
  zz_1(gg_1)=\tau\tau_1,
  \]
  where $zz_1\in C$, $gg_1\in \Gamma_1$ and $\tau\tau_1\in \Gamma_\alpha$. Since $zz_1$ is not the identity, the argument above
  shows that each element $zz_1$, $gg_1$ and $\tau\tau_1$ has order 2. But $g$ and $g_1$ are elements of order 2 in a dihedral group
  of order $2n$, where $n$ is odd, and it is well known and easy to prove that their product has odd order. We have arrived at a contradiction
  and we have thus proved that $|C|=2$, as required. 
  
  To finish the argument, $\Gamma/C$ is isomorphic to a subgroup of the automorphism group of $\Gamma_1$ and a comparison of orders, taken
  in conjunction with the conclusion of Lemma \ref{automorphism_group_of_dihedral}, shows that $\Gamma/C$ is isomorphic to the entire
  automorphism group of $\Gamma_1$. $C$ is central in $\Gamma$, since it is normal and has order 2. We are required to show that it is a direct factor. Let $z$ be a generator of $C$. As we have already argued, $z$ fixes no roots of $D(T)$, since it centralizes $\Gamma_1$. Thus
  in its permutation action on the $2n$ roots, $z$ is a product of $n$ disjoint transpositions. Since $n$ is odd, $z$ determines an odd permutation
  and hence in the symmetric group of degree $2n$, $\Gamma$ is the direct product of a subgroup of index 2, lying in the alternating group
  of degree $2n$, and the central subgroup generated by  the permutation corresponding to $z$. This completes the proof.
\end{proof}

\section{Further Descent Analysis in  Characteristic $p$}\label{ff}

In this section, 
we will realize the automorphism groups of several subgroups of $PGL(2,q)$ as Galois groups
over $\mathbb{F}_p(x)$.
We do this by applying the descent methods outlined in the 
previous section to go from $\F_q(x)$ to $\F_p(x)$. 
We start by realizing the automorphism group of $PGL(2,q)$ as a Galois group over $\mathbb{F}_p(x)$.
This group is the same as the automorphism group of $PSL(2,q)$. 
We also realize a subgroup of index 2 in the automorphism group
of  $PSL(2,q)$ as a Galois group over $\mathbb{F}_p(x)$.
The automorphism group of $PGL(2,q)$ is sometimes known as $P\Gamma L(2,q)$,
as discussed in \cite{Ab}, Section 15, for example.

We briefly recall some properties of $PGL(2,q)$ that we need to prove subsequent results. If $q$ is odd and at least 5, $PGL(2,q)$ has a unique proper normal subgroup $PSL(2,q)$, which has index 2 in $PGL(2,q)$. $PSL(2,q)$ is a nonabelian simple group. The centre of $PGL(2,q)$ is  trivial. When $q=3$, $PGL(2,3)$ is isomorphic to the symmetric group $S_4$, and $PSL(2,3)$ is isomorphic to the alternating group
$A_4$, whose centre is also trivial. $PSL(2,3)$ contains a normal abelian subgroup of order 4, which is  normal in $PGL(2,3)$.

When $q$ is a power of 2, $PGL(2,q)$ is a nonabelian simple group provided that $q>2$. 

\begin{lemma} \label{subgroup_of_automorphism_group}

Let $q=p^n$, where $p$ is a prime, and let $X=Aut(PGL(2,q))$ denote the automorphism group of $PGL(2,q)$. Likewise, let $Y=Aut(PSL(2,q))$ denote the automorphism group of $PSL(2,q)$. Then $X$ is naturally a subgroup of $Y$ and hence has order dividing $nq(q^2-1)$.
\end{lemma}

\begin{proof}
We set $G=PGL(2,q)$ and $H=PSL(2,q)$. We begin by considering the case where $q$ is a power of 2. In this circumstance, $G=H$, and thus
automatically $X=Y$. We may therefore make the assumption that $q$ is odd.

It is well known that $H$ is the commutator subgroup of $G$ and hence is characteristic in $G$. Thus, restriction of the elements of 
$X$ to $H$ naturally defines a homomorphism from $X$ into $Y$. We claim that if an element $\sigma\in X$
acts trivially on $H$, then $\sigma$ is the identity.

To prove the claim, we proceed as follows. Suppose that $\sigma(h)=h$ for all $h$ in $H$. Let $t$ be an element of $G\setminus H$.
We want to show that $\sigma(t)=t$ also, from which our claim will follow. Since $t^{-1}ht\in H$, we have
\[
\sigma(t^{-1}ht)=t^{-1}ht=\sigma(t)^{-1}\sigma(h)\sigma(t)=\sigma(t)^{-1}h\sigma(t).
\]
This shows that $\sigma(t)t^{-1}$ centralizes $H$ and hence is an element of
the centralizer of $H$ in $G$, which we shall call $C$.

It is certainly the case that $C$ is normal in $G$, since $H$ is normal in $G$. Now, as it is well known that $H$ is the unique proper normal subgroup of $G$, except when $q=3$, either $C=1$ or
$C=H$ when $q>3$. 
However, $H$ is nonabelian and simple when
$q>3$, and hence we cannot have $C=H$, because this implies that $H$ is abelian. It follows that $C=1$ and so $\sigma(t)=t$. This establishes the claim for $q>3$.

We need to examine what happens when $q=3$. It is still the case that $H$ contains all proper normal subgroups of $G$ and thus $C$ lies in $H$. Then $C$ must  be the centre of $H$. But when $q=3$, as we noted above, $H$ is isomorphic to $A_4$ and thus has trivial centre. This forces the conclusion that $C=1$ and we again have that $\sigma$ is the identity.

Having established this embedding of $X$ into $Y$, our estimate for the order of $X$ follows for example from Theorem 3.2 of \cite{W}. 
\end{proof}

We remark that we could be more precise in the previous lemma, since $X$ equals $Y$, but this will become apparent in the proof that follows.

\subsection{$PGL(2,q)$}

We now apply this analysis to identify the Galois group $\Gamma$ in Lemma \ref{split} when $G=PGL(2,q)$. Before embarking on the proof, we
extract a lemma needed in the course of the argument.

\begin{lemma} \label{order_of_centralizer}
Let $G$ be a subgroup of $PGL(2,q)$, where $q$ is a power of the prime $p$. Suppose that the centralizer in $G$ of some nonidentity
element $g$ has order divisible by $p$. Then $g$ has order $p$ and the order of the centralizer is a power of $p$.
\end{lemma}

This result must be well known, and it is easy to prove. For a formal proof, we refer to Lemma 4 of \cite{GM}.

  \begin{thm} \label{automorphism_group}
  Let $q=p^n$, where $p$ is a prime, and let $G=PGL(2,q)$. Let $\Phi(x)=f(x)/g(x)$ be a nonconstant coefficient in the orbit polynomial
  of $G$. Let $P=f(T)-xg(T)$ and consider $P$ as an element of $\mathbb{F}_p(x)[T]$. Then 
  the Galois group $\Gamma$ of $P$ over $\mathbb{F}_p(x)$ is isomorphic to the automorphism group $Aut(G)$ of $G$ and has order $n|G|$.
  \end{thm}
  
  \begin{proof}
  
  We employ the notation of Hypotheses (H), with $F=\mathbb{F}_q$, $E=\mathbb{F}_p$, $G=PGL(2,q)$ and $H$ cyclic of order $n$. Let $\sigma$ be a nonidentity element of $\Gamma_\alpha$. We have seen that $\sigma$ induces a field automorphism of $\mathbb{F}_q$. Let $p^a$ be the order of the subfield of $\mathbb{F}_q$ fixed elementwise by $\sigma$. Now $\sigma$ also induces a group automorphism of $PGL(2,q)$ and it is straightforward
  to see that the subgroup of fixed points of this action is $PGL(2,p^a)$. We note here that the order of this subgroup
  is divisible by $p$. Lemma \ref{number_of_roots_fixed} tells us that the number of roots of $P$ fixed by $\sigma$ is $|PGL(2,p^a)|$.
  
  Turning to the Galois group $\Gamma$, we know that $\Gamma$ contains a normal subgroup $\Gamma_1$ isomorphic to $PGL(2,q)$, and the quotient
  is cyclic of order $n$. Let $C$ denote the centralizer of $\Gamma_1$ in $\Gamma$. Since the centre of $PGL(2,q)$ is trivial, $\Gamma_1\cap C=1$. Suppose that $C$ is nontrivial. Then by Lemma \ref{pinning_down_the_centralizer}, there is a nonidentity element $\tau$, say, of
  $\Gamma_\alpha$ such that the number of roots of $P$ fixed by $\tau$ is the order of the centralizer in $\Gamma_1$ of some nonidentity
  element $g$ of $\Gamma_1$. 
  
  The argument above shows that the order of the centralizer of $g$ in $\Gamma_1$ is divisible by $p$. Since $\Gamma_1$ is isomorphic to
  $PGL(2,q)$, Lemma \ref{order_of_centralizer} tells us that $g$ has order $p$ and the order of the centralizer is a power of $p$.
  This in turn implies that $|PGL(2,p^a)|$ is a power of $p$, which is certainly not true. It follows that $C$ is trivial.

  As is well known, since $\Gamma_1$ is normal in $\Gamma$, we have a homomorphism from $\Gamma$ into the automorphism group of $\Gamma_1$, defined by conjugation action, and the kernel
  is $C$. Given the triviality of $C$ just proved, $\Gamma$ is realized as a subgroup
  of $Aut(\Gamma_1)$. We also know from Lemma \ref{subgroup_of_automorphism_group} that $|Aut(\Gamma_1)|$ has order dividing $n|G|$. Since $|\Gamma|=n|\Gamma_1|$, it follows that $\Gamma$ is isomorphic to $Aut(PGL(2,q))$, as claimed.
  \end{proof}
  
  \subsection{$PSL(2,q)$}

  There is an analogue of Theorem \ref{automorphism_group} that applies to the group $PSL(2,q)$. As we mentioned earlier,
  we realize not the whole automorphism group but rather a subgroup of index 2. We confine attention to the case that $q$ is odd, since the groups $PGL(2,q)$ and $PSL(2,q)$ are identical when $q$ is a power of 2.
  
  \begin{thm} \label{small_automorphism_group}
  Let $q=p^n$, where $p$ is an odd prime, and let $G=PSL(2,q)$. Let $\Phi(x)=f(x)/g(x)$ be a nonconstant coefficient in the orbit polynomial
  of $G$. Let $C(T)=f(T)-xg(T)$ and consider $C(T)$ as an element of $\mathbb{F}_p(x)[T]$. Then 
  the Galois group $\Gamma$ of $C(T)$ over $\mathbb{F}_p(x)$ is isomorphic to a subgroup of the automorphism group  of $G$ and has order $n|G|$.
  \end{thm}
  
  \begin{proof}
  The proof is essentially identical to the proof of the previous theorem. The main point of difference is the following. Let
  $\alpha$ be a root of $C(T)$ in a splitting field over $\mathbb{F}_p(x)$. Then the other roots of $C(T)$ have the form
  \[
  \frac{a\alpha+b}{c\alpha+d},
  \]
  where $a$, \dots, $d$ are elements of $\mathbb{F}_q$ with $ad-bc$ a nonzero square in $\mathbb{F}_q$. When an element $\tau$ of $\Gamma_\alpha$ acts on the roots, the number of roots it fixes is either $|PSL(2,p^a)|$ or $|PGL(2,p^a)|$, where $p^a$ is the order of the subfield of fixed points of $\tau$ acting on $\mathbb{F}_q$. The vital part of the argument is that the number of fixed roots is still divisible by $p$
  and this enables the rest of the proof to proceed as in the former case.
  \end{proof}
  
  \subsection{Affine General Linear Group}
  
  Similar techniques can be applied to realize another  subgroup of $PGL(2,q)$ as a Galois group over $\mathbb{F}_q(x)$ and to realize its full automorphism group as a Galois group over 
  $\mathbb{F}_p(x)$. We assume as usual that $q=p^n$, where $p$ is a prime.
  
  Consider the subgroup $G$, say, of $PGL(2,q)$ consisting of all affine linear transformations of the form $x\mapsto ax+b$, where $a$ and $b$ run over all elements in $\mathbb{F}_q$, subject to $a\neq 0$. It is clear that $|G|=q(q-1)$. $G$ is usually called the (one-dimensional) affine
  general linear group and is denoted by $AGL_1(\mathbb{F}_q)$. $G$ is the normalizer of a Sylow $p$-subgroup of $PGL(2,q)$.
  
  We note that $G$ is invariant under the action of the Galois group of $\mathbb{F}_q$ over $\mathbb{F}_p$. Thus we have an action of the cyclic group of order $n$ as automorphisms of $G$. Aside from the case when $q=2$, $G$ has trivial centre and thus $G$ coincides with its subgroup of inner automorphisms. The Galois group induces a cyclic subgroup of automorphisms of $G$, which, apart from the identity, are not inner.
  Thus $G$ has a group of automorphisms of order $nq(q-1)$, containing a normal subgroup isomorphic to $G$, and having a cyclic quotient of order $n$. 
  The group of automorphisms obtained is denoted by $A\Gamma L_1(\mathbb{F}_q)$ and it can be 
  shown\footnote{For a proof, provided by Derek Holt, we refer to the question \lq\lq Automorphism group of the general affine group of the affine line over a finite field\rq\rq, posted on StackExchange and answered on 19 February 2015.}
  to be the full group of automorphisms of $AGL_1(\mathbb{F}_q)$, except when $q=2$.
    
  When $G$ acts on $\mathbb{F}_q(x)$, it leaves invariant the polynomial $(x^q-x)^{q-1}$, as is easily confirmed. Since this polynomial
  has degree $q(q-1)=|G|$, it generates the field of $G$-invariant rational functions. It follows from Theorem \ref{main1} that the polynomial
  $(T^q-T)^{q-1}-x$ in $\mathbb{F}_q(x)[T]$ has Galois group over $\mathbb{F}_q(x)$ isomorphic to $G$. 
  
  We show that over $\mathbb{F}_p(x)$, $(T^q-T)^{q-1}-x$ has Galois group isomorphic to $A\Gamma L_1(\mathbb{F}_q)$ and we thus realize
  the automorphism group of $AGL_1(\mathbb{F}_q)$ as a Galois group over $\mathbb{F}_p(x)$, the case $q=2$ being excluded.
  
  \begin{thm} \label{affine_group}
  Let $q=p^n$, where $p$ is a prime and $q>2$. Let $B(T)=(T^q-T)^{q-1}-x$ and consider $B(T)$ as an element of $\mathbb{F}_p(x)[T]$. Then  the Galois group $\Gamma$ of $B(T)$ over $\mathbb{F}_p(x)$ is isomorphic to the automorphism group of the affine general linear group
  $AGL_1(\mathbb{F}_q)$ and has order $nq(q-1)$.
  \end{thm}
  
  \begin{proof}
  We follow the proof of Theorem \ref{automorphism_group} and only give a sketch, as the proofs are almost identical. We again adopt
  the notation of Hypotheses (H), with $F=\mathbb{F}_q$, $E=\mathbb{F}_p$, $G=AGL_1(\mathbb{F}_q)$ and $H$ cyclic of order $n$.
  $\Gamma$ contains a normal subgroup $\Gamma_1$  isomorphic to $AGL_1(\mathbb{F}_q)$, with quotient cyclic of order $n$.

  Let $\alpha$ be a root of $B(T)$ in the splitting field.  The roots of $B(T)$ have the form
  $a\alpha+b$, where $a$ and $b$ run over $\mathbb{F}_q$, with $a\neq 0$. Let $\sigma$ be an element of $\Gamma_\alpha$ and let $\beta=
  a\alpha+b$ be a root of $B(T)$. Then
  \[
  \sigma(\beta)=\sigma(a)\alpha+\sigma(b)
  \]
  and it follows that the number of roots fixed by $\sigma$ is $p^a(p^a-1)$, where $p^a$ is the order of the subfield of $\mathbb{F}_q$
  fixed elementwise by $\sigma$. 
  
  Let $C$ be the centralizer of $\Gamma_1$ in $\Gamma$. Since $\Gamma_1$ has trivial centre when $q>2$, $\Gamma_1\cap C=1$. Suppose that $C$ is nontrivial. Then by Lemma \ref{pinning_down_the_centralizer}, there is a nonidentity element $\tau$, say, of
  $\Gamma_\alpha$ such that the number of roots of $B(T)$ fixed by $\tau$ is the order of the centralizer in $\Gamma_1$ of some nonidentity
  element $g$ of $\Gamma_1$. 
  
  The argument above now shows that the order of the centralizer of $g$ in $\Gamma_1$ has the form 
  $p^a(p^a-1)$ and in turn
   Lemma \ref{order_of_centralizer} tells us that $g$ has order $p$ and the order of the centralizer is a power of $p$. Since $p^a(p^a-1)$
  is not a power of $p$, we have a contradiction and we have therefore shown that $C=1$.

  Consequently, $\Gamma$ embeds into the automorphism group of 
  $AGL_1(\mathbb{F}_q)$ and since the automorphism group
  has the same order as $\Gamma$, namely $nq(q-1)$,
   $\Gamma$ is isomorphic with the automorphism group of $AGL_1(\mathbb{F}_q)$ .
  \end{proof}

\section{Further Remarks}

Apart from the descent idea, there are other methods to reach groups that are not 
subgroups of $PGL(2,F)$, building on
a subgroup that has been realized. Here is one approach.

Given a subgroup $G$ of $PGL(2,F)$ we have a generator $\Phi=f/g$ of the field of $G$-invariant functions. In Theorem \ref{main1}  we proved that $f(T)-xg(T)$ has Galois group $G$ over $F(x)$.
We also showed in Lemma 
that if $f$ and $g$ are relatively prime, then $f(T)-xg(T)$ is irreducible.
If $f$ and $g$ are relatively prime, so are $f(T^2)$ and $g(T^2)$, so 
$f(T^2)-xg(T^2)$ is irreducible as well.

Let $\Gamma$ be the Galois group of $f(T^2)-xg(T^2)$. 
There is an elementary abelian normal 2-subgroup $A$ of $\Gamma$ such that $\Gamma/A\cong G$.
As the polynomial $f(T^2)-xg(T^2)$ has degree $2|G|$, $2|G|$ divides $|\Gamma|$.
This means that $A$ is nontrivial.

Example: let $G=PSL(2,3)\cong A_4$, we calculate using a formula given by Carlitz that
$f(x)=x^{12} - x^{10} + x^6 - x^2 + 1$ and $g(x)=x^9 -x^3$, where $\Phi(x)=f(x)/g(x)$.
Then we found using Magma that the Galois group of 
$f(T^2)-xg(T^2)=T^{24} - T^{20} + T^{12} - T^4 + 1 +x(T^{18} -T^6)$
is a nonsplit extension of $C_2^4$ by $A_4$ of order 192.
In this case the subgroup $A$ has order 16.

\end{document}